\documentclass[11pt]{article}

\usepackage{subfigure}

\usepackage{graphicx}
\usepackage{amsmath, amsthm, amssymb}
\usepackage{graphicx}
\input epsf.tex
\usepackage{epsf}
\usepackage{epsfig}
\usepackage[vcentering,dvips]{geometry}
\newtheorem{theorem}{Theorem}[section]
\newtheorem{lemma}[theorem]{Lemma}
\newtheorem{prop}[theorem]{Proposition}
\newtheorem{corollary}[theorem]{Corollary}
\newtheorem{definition}[theorem]{Definition}

\newtheorem{remark}{Remark}

\newcommand{\sgn}{\operatorname{sgn}}

\usepackage{fullpage}
\def\nnew{\color{red}}
\def\mnew{\color{black}}
\def\nnew{}
\def\mnew{}

\begin{document}

% paper title
\title{Iterative Thresholding meets Free Discontinuity Problems}
\author{
Massimo Fornasier\footnote{Johann Radon Institute for Computational and 
Applied Mathematics, Austrian Academy of Sciences, Altenbergerstrasse 69, A-4040 Linz, Austria,
email: {\tt massimo.fornasier@oeaw.ac.at}.} \ and  
Rachel Ward\footnote{Program in Applied and Computational Mathematics, Princeton University, Fine Hall, Washington Road, 08544 Princeton, NJ, U.S.A., email:  
{\tt rward@math.princeton.edu}.}}

% make the title area
\maketitle

\begin{abstract}
Free-discontinuity problems describe situations where the solution of interest is defined by a function and a lower dimensional set consisting of the discontinuities of the function.  Hence, the derivative of the solution is assumed to be a `small' function almost everywhere except on sets where it concentrates as a singular measure. This is the case, for instance, in crack detection from fracture mechanics or in certain digital image segmentation problems. If we discretize such situations for numerical purposes, the free-discontinuity problem in the discrete setting can be re-formulated as that of finding a derivative vector with small components at all but a few entries that exceed a certain threshold.  This problem is similar to those encountered in the field of `sparse recovery', where vectors with a small number of dominating components in absolute value are recovered from a few given linear measurements via the minimization of related energy functionals. Several iterat!
 ive thresholding algorithms that intertwine gradient-type iterations with thresholding steps have been designed to recover sparse solutions in this setting.  It is natural to wonder if and/or how such algorithms can be used towards solving discrete free-discontinuity problems.  The current paper explores this connection, and, by establishing an iterative thresholding algorithm for discrete free-discontinuity problems, provides new insights on properties of minimizing solutions thereof.
\end{abstract}

\\


\tableofcontents
\section{Introduction}

In the following introductory sections, we will establish the mathematical setting of the paper, and review the features of free-discontinuity problems that are relevant to the current discussion.

\subsection{Free-discontinuity problems: the Mumford-Shah functional}
The terminology `free-discontinuity problem' was introduced by De Giorgi \cite{dg91}
%#\footnote{E. De Giorgi, \emph{ Free-discontinuity problems in calculus of variations}, Frontiers in pure and applied mathematics, a collection of papers dedicated to J.-L. Lions on the occasion of his $60^{th}$ birthday, R. Dautray, ed. North Holland, 1991, 55-62.} 
to indicate a class of variational problems that consist in the minimization of a functional, involving both volume and surface energies, depending on a closed set $K \subset \mathbb R^d$, and a function $u$ on $\mathbb R^d$ usually smooth outside of $K$. In particular, 
\begin{itemize}
\item $K$ is not fixed a priori and is an unknown of the problem;
\item $K$ is not a boundary in general, but a free-surface inside the domain of the problem.
\end{itemize}

\noindent The best-known example of a free-discontinuity problem is the one modelled by the so-called Mumford-Shah functional \cite{MS89}, which is defined by
$$
J(u,K):= \int_{\Omega \setminus K} \left [ | \nabla u |^2 + \alpha ( u - g )^2 \right ] dx +\beta \mathcal{H}^{d-1}(K \cap \Omega).
$$
The set $\Omega$ is a bounded open subset of $\mathbb{R}^d$, $\alpha, \beta >0$ are fixed constants, and $g \in L^\infty(\Omega)$. Here $\mathcal H^N$ denotes the $N$-dimensional Hausdorff measure. Throughout this paper, the dimension of the underlying Euclidean space $\mathbb R^d$ will always be $d=1$ or $d=2$.  In the context of visual analysis, $g$ is a given noisy image that we want to approximate by the minimizing function $u \in W^{1,2}(\Omega \setminus K)$; the set $K$ is simultaneously used in order to {\it segment} the image into connected components.  For a broad overview on free-discontinuity problems, their analysis, and applications, we refer the reader to \cite{AFP}.
\\

\noindent If the set $K$ were fixed, then the minimization of $J$ with respect to $u$ would be a relatively simple problem, equivalent to solving the following system of equations:
\begin{eqnarray*}
\Delta u&=& \alpha ( u - g ), \qquad \mbox{in } \Omega \setminus K,\\
\frac{\partial u}{\partial \nu} &=&0, \qquad  \qquad \quad \mbox{ on } \partial \Omega \cup K,   
\end{eqnarray*}
where $\nu$ is the outward-pointing normal vector at any $x \in \partial \Omega \cup K$. Therefore the relevant unknown in free-discontinuity problems is the set $K$.
Ensuring the existence of minimizers $(u,K)$ of $J$ is a challenging problem because there is no topology on the closed sets that ensures 
\begin{itemize} 
\item[(a)] compactness of minimizing sequences and 
\item[(b)] lower semicontinuity of the Hausdorff measure. 
\end{itemize}
Indeed, it is well-known, by the direct method of calculus of variations \cite[Chapter 1]{DM}, that the two previous conditions ensure the existence of minimizers.  However, the problem becomes more manageable if we restrict our domain to functions $u \in BV(\Omega) \cap  W^{1,2}(\Omega \setminus K)$, and make the identification $K \equiv \overline{S_u}$ where $S_u$ is the well-defined discontinuity set of $u$.  In this case, we need to work only with a topology on the space $BV (\Omega)$ of bounded variation, and no set topology is anymore required. \\ \\
Unfortunately the space $BV(\Omega)$ is `too large'; it contains Cantor-like functions whose approximate gradient vanishes, $\nabla u= 0$, almost everywhere, and whose discontinuity set has measure zero, $\mathcal{H}^{d-1}(S_u)=0$.  As these functions are dense in $L^2(\Omega)$, the problem is trivialized; see \cite{AFP} for details.
\\
\\
Nevertheless, it is possible to give a meaningful formulation of the functional $J$ if we exclude such functions and restrict $J$ to the space $SBV(\Omega)$ constituted of $BV$-functions with vanishing Cantor part.   If we assume again $K \equiv \overline{S_u}$,  the solution can be recast as the minimization of
\begin{equation}
\label{amb}
\mathcal J ( u ) = \int_{\Omega \setminus S_u} \left [ | \nabla u |^2 + \alpha ( u - g )^2 \right ] dx +\beta \mathcal{H}^{d-1}(S_u).
\end{equation}
The existence of minimizers in $SBV$ for the functional \eqref{amb} was established by Ambrosio on the basis of his fundamental compactness theorem in \cite{AM89}, see also \cite[Theorem 4.7 and Theorem 4.8]{AFP}.

\subsection{$\Gamma$-convergence approximation to free-discontinuity problems}

The discontinuity set $S_u$ of a $SBV$-function $u$ is not an object that can be easily handled, especially numerically. This difficulty gave rise to the development of approximation methods for the Mumford-Shah functional and its minimizers where \emph{sets} are no longer involved, and instead substituted by suitable \emph{indicator functions}.  In order to understand the theoretical basis for these approximations, we need to introduce the notion of $\Gamma$-convergence, which is today considered one of the most successful notions of `variational convergence'; we state only the definition of $\Gamma$-convergence below, but refer the reader to \cite{DM,br02} for a broad introduction.
\begin{definition}
Let $(X,d)$ be a metric space\footnote{ Observe that by \cite[Proposition 8.7]{DM} suitable bounded sets $X$ endowed with the weak topology induced by a larger Banach space are indeed metrizable, so this condition is not that restrictive.} and let $f,f_n : X \to [0,\infty]$ be functions for $n \in \mathbb N$. We say that $(f_n)_{n \in \mathbb N}$ $\Gamma$-converges to $f$ if the following two conditions are satisfied:
\begin{itemize}
\item[i)] for any sequence $(x_n)_n \subset X$ converging to $x$, 
$$
\liminf_n f_n(x_n) \geq f(x);
$$
\item[ii)] for any $x \in X$, there exists a sequence $(x_n)_n \subset X$ converging to $x$ such that
$$
\limsup_n f_n(x_n) \leq f(x).
$$

\end{itemize}
\label{GammaDef}
\end{definition}

One important consequence of Definition \ref{GammaDef} is that if a sequence of functionals $f_n$ $\Gamma$-converges to a target functional $f$, then the corresponding minimizers of $f_n$ also converge to minimizers of $f$, see \cite[Corollary 7.30]{DM}. 
\\
\\
We define now
\begin{equation}
F_{\varepsilon}(u,v) := \int_{\Omega} \left [ v^2 | \nabla u |^2 + \alpha ( u - g )^2 \right ]  +\frac{\beta}{2 }  \left ( \varepsilon |\nabla v|^2 + \frac{(1 -v)^2}{\varepsilon} \right ) dx
\end{equation}
over the domain $(L^2(\Omega))^2$, along with the related functional
\begin{equation}
\label{amto}
\mathcal J_\varepsilon(u,v) :=   \left\{\begin{array}{cl} F_\varepsilon(u,v)
	  & \textrm{, if } v \in W^{1,2}(\Omega) \textrm{, }u v \in W^{1,2}(\Omega) \textrm{, and } 0 \leq v \leq 1,\\
	\infty   & \textrm{, else.}
	   \end{array}\right.
\end{equation}
Note that at the minimizer $(u,v)$ of $\mathcal J_\varepsilon$, the function $0 \leq v \leq 1$ tends to indicate the discontinuity set $S_u$ of the functional \eqref{amb} as $\varepsilon \to 0$. In \cite{amto90} Ambrosio and Tortorelli proved the following $\Gamma$-approximation result:
\begin{theorem}[Ambrosio-Tortorelli '90] 
For any infinitesimal sequence $(\varepsilon_n)_n$, the functional $\mathcal J_{\varepsilon_n}(u,v)$ $\Gamma$-converges in $(L^2(\Omega))^2$ to the functional
\begin{equation}
\mathcal J ( u,v):=  \left\{\begin{array}{cl}  \mathcal J(u) & \textrm{, if } v \equiv 1, \\
 \infty & \textrm{, otherwise.}
 \end{array}\right.
 \end{equation}
\end{theorem}

\subsection{Discrete approximation}

In fact, the Mumford-Shah functional is the continuous version of a previous discrete formulation of the image segmentation problem proposed by Geman and Geman in \cite{GG}; see also the work of Blake and Zisserman in \cite{BZ}. Let us recall this discrete approach.
For simplicity let $d=2$ (as for image processing problems),  $\Omega = [0,1]^2$, and let $u_{i,j}=u(h i,h j)$, $(i,j) \in \mathbb{Z}^2$ be a discrete function defined on $\Omega_h:=\Omega \cap h \mathbb{Z}^2$, for $h>0$.
Define $W_h ( t ) = \min\{ t^2, \beta/h\}$ to be the truncated quadratic potential, and
\begin{eqnarray*}
\mathcal J_{\sqrt{\beta/h}}(u ) &:=& h^2 \sum_{(h i, h j) \in \Omega_h} W_h \left ( \frac{u_{i+1,j} - u_{i,j}}{h} \right ) \\
&+& h^2 \sum_{ ( h i, h j) \in \Omega_h} W_h \left ( \frac{u_{i,j+1} - u_{i,j}}{h} \right )\\
&+& \alpha h^2 \sum_{(h i,h j) \in \Omega_h} ( u_{i,j} - g_{i,j})^2.
\label{discrete}
\end{eqnarray*}
Chambolle \cite{ca92,ca95} gave formal clarification as to how the discrete functional $\mathcal J_{\sqrt{\beta/h}}$ approximates the continuous functional $\mathcal J$ of Ambrosio: discrete sequences can be interpolated by piecewise linear functions in such a way as to allow for discontinuities when the discrete finite differences of the sampling values are large enough. On the basis of this identification of discrete functions on $\Omega_h$ and functions defined on the `continuous domain' $\Omega$, we have the following result:
\begin{theorem}[Chambolle '95]
\label{chthm} 
The functional $\mathcal J_{\sqrt{\beta/h}}$ $\Gamma$-converges in $\mathcal B(\Omega) $ (the space of Borel-measurable functions, which is metrizable, see \cite{ca95} for details) to 
$$
\mathcal J^{cab} ( u ) = \int_{\Omega \setminus S_u} \left [ | \nabla u |^2 + \alpha ( u - g )^2 \right ] dx +\beta \mathcal{C}(S_u),
$$
as $h \to 0$, where $\mathcal{C}$ is the so-called `cab-driver' measure defined below.
\end{theorem}
Basically $\mathcal{C}$ measures the length of a curve only through its projections along horizontal and vertical axes; for a regular $C^1$ curve $c=\gamma([0,1])$, with $\gamma(t)=(\gamma_1(t), \gamma_2(t)) \in \Omega$, we have
$$
\mathcal{C}(c) = \int_0^1 \left ( |\gamma_1'(t)| + |\gamma_2'(t)| \right) dt.
$$
The reason this anisotropic (or, direction dependent) measure appears, in place of the Hausdorff measure in the Mumford-Shah functional, is due to the approximation of derivatives by finite differences defined on a `rigid' squared geometry. A discretization of derivatives based on meshes {\it adapted} to the morphology of the discontinuity indeed leads to precise approximations of the Mumford-Shah functional  \cite{chdm99,boch00}.  

\subsection{Free-discontinuity problems and discrete derivatives}
 In the literature, several methods have been proposed to numerically approximate minimizers of the Mumford-Shah functional \cite{beco94,boch00,ca92,ca95, ma92}.
In particular, a relaxation algorithm, based essentially on alternated minimization of a finite element approximation of the Ambrosio and Tortorelli functional \eqref{amto}, leads to iterated solutions of suitable elliptic PDEs, where the differential part includes the auxiliary variable $v$ which encodes and indicates information about the discontinuity set. These implementations are basically finite dimensional approximations to the following algorithm: \emph{Starting with $v^{(0)} \equiv 1$, iterate}
$$
\left \{ \begin{array}{l}
u^{(n+1)} := \arg \min_{u \in W^{1,2}(\Omega)} \mathcal J_\varepsilon(u ,v^{(n)})\\
v^{(n+1)} := \arg \min_{v \in W^{1,2}(\Omega)} \mathcal J_\varepsilon(u^{(n+1)} ,v).
\end{array}
\right .
$$
{  However, neither has a proof of convergence of this iterative process to its stationary points been explicitly provided in the literature, nor have the properties of such stationary points been investigated, especially in case of genuine inverse problems (see the discussion in Subsection \ref{MS4invprob}).

In this paper, we take a different approach and investigate how minimization of the $\Gamma$-approximating discrete functionals \eqref{discrete} can be implemented efficiently by iterative thresholding on the discrete derivatives. 
Unlike the aforementioned approach, we will be able to provide a rigorous proof of convergence to stationary points, which coincide with local minimizers of the discrete Mumford-Shah functional.  Moreover, we are able to characterize stability properties of such stationary points, and demonstrate the stability of global minimizers of the discrete Mumford Shah functional.}

Let us recall: the solutions $u$ of a free-discontinuity problem are supposed to be smooth out of a minimal ipersurface $K$. This means that the distributional derivative of $u$ is a `small function' everywhere except on $K$ where it coincides with a singular measure. In the discrete approximation $\eqref{discrete}$, the vector of finite differences $(w_j) = (\frac{u_{i,j+1} - u_{i,j}}{h}, \frac{u_{i+1,j} - u_{i,j}}{h} )$ corresponds to a piecewise constant function that is small everywhere except for a few locations, corresponding to $|w_j| \geq \sqrt{\beta/h}$, that approximate the discontinuity set $K$.   So, in terms of derivatives, solutions of $\eqref{discrete}$ are vectors having only few large entries.  In the next section, we clarify how we can indeed work with just derivatives and forget the primal problem.
\subsubsection{The 1-D case}
Let us assume for simplicity that the dimension $d=1$, the domain $\Omega = [0,1]$, and the parameters $\alpha=\beta=1$.  Denote by $u_{i}=u(h i)$ a discrete function defined on $h i \in \Omega_h:=\Omega \cap h \mathbb{Z}$, for $h>0$; note that the vector $(u_i) \in \mathbb{R}^n$ for $n = \lfloor1/h \rfloor$. In this setting, the discrete functional $\eqref{discrete}$ reduces to
\begin{eqnarray*}
\mathcal J_{\sqrt{1/h}}(u ) &=& h \sum_{(h i) \in \Omega_h} W_h \left ( \frac{u_{i+1} - u_{i}}{h} \right ) \\
&+& h \sum_{(h i) \in \Omega_h} ( u_{i} - g_{i})^2,
\end{eqnarray*}
where we recall that $W_h(t) = \min\{t^2, 1/h \}$.  Since no geometrical anisotropy is now involved ($d=1$), it is possible to show that this discrete functional $\Gamma$-converges precisely  to the corresponding Mumford-Shah functional on intervals \cite{ca92}.
\\
\\
For $(u_i)_{ hi \in \Omega_h}$ we define the discrete derivative as the matrix $D_h: \mathbb{R}^n \rightarrow \mathbb{R}^{n-1}$ that maps $(u_i)_{ hi \in \Omega_h}$ into $\left (\frac{u_{i+1}-u_i}{h} \right )_{i}$, given by
\begin{equation}
\label{dermtrx}
D_h = \frac{1}{h} \left ( \begin{array}{cccccc} -1 & 1 & 0& \dots & \dots &0\\
0&-1&1&0 &\dots&0\\
\dots\\
0&0&\dots&\dots&-1&1
\end{array} \right ).
\end{equation}
It is not too difficult to show that
$$
 u =D^\dagger_h D_h u + c,
$$
where $D_h^\dagger$ is the pseudo-inverse matrix of $D_h$ (in the Moore-Penrose sense; note that $D^\dagger_h$ maps $\mathbb R^{n-1}$ into $\mathbb R^n$ and is an injective operator) and $c$ is a constant vector which depends on $u$, and the values of its entries coincide with the mean value $h \sum_{ hi \in \Omega_h} u_i$ of $u$.
Therefore, any vector $u$ is uniquely identified by the pair $(D_h u , c)$.
\\

 \noindent Since constant vectors comprise the null space of $D_h$, the orthogonality relation $\langle D^\dagger_h D_h u,c\rangle_{\ell_2^n} =0$ holds for any vector $u$ and any constant  vector $c$. Here the scalar product $\langle \cdot, \cdot \rangle_{\ell_2^n} = \sum_i u_i v_i$ is the standard Euclidean scalar product, which induces the Euclidean norm $\| u \|_{\ell_2^n} := \left ( \sum_ i u_i^2 \right )^{1/2}$. 
Using this orthogonality property, we have that
\begin{eqnarray*}
\| u - g \|_{\ell_2^n}^2 &=& \| D^\dagger_h D_h u -  D^\dagger_h D_h g + ( c - c_g )\|_{\ell_2^n}^2 \\
&=&    \| D^\dagger_h D_h u -  D^\dagger_h D_h g\|_{\ell_2^n}^2 + \| c - c_g \|_{\ell_2^n}^2 
\end{eqnarray*}
Hence, with a slight abuse of notation, we can reformulate the original problem in terms of derivatives, and mean values, by
\begin{eqnarray*}
\mathcal J_{1/\sqrt h}(z,c) &=& h \| D_h^\dagger z  - f \|_{\ell_2^{n}}^2 + h \| c- c_{g} \|_{\ell_2^n}^2 + h \sum_{i} \min \left \{|z_i|^2, \frac{1}{h}  \right \} \\
\end{eqnarray*}
where $z = D_h u$ and $f =  D^\dagger_h D_h g$. Of course at the minimizer $u$ we have $c = c_g$, since this term in $\mathcal J_{1/\sqrt h}$ does not depend on $z$. Therefore, $\| c - c_g \|_2^2$ does not play any role in the minimization and can be neglected. Once the minimal derivative vector $z$ is computed, we can assemble the minimal $u$ by incorporating the mean value of $g$ as follows:
$$
 u = D_h^\dagger z + c_g.
$$
%\begin{remark}[$\circledast \circledast$]
% In the following we will introduce and study iterative thresholding algorithms which minimize functionals of the type presented above. The most onerous effort in the iterations is due to the application of the bounded operator which is involved in the discrepancy with respect to data in the functional. In the case described above, this operator is $D^\dagger_h$.
% So, for the algorithms to be efficient or even just implementable, it is necessary to have matrices/operators which are easily stored in memory and whose application is not too expensive computationally. On the one hand, for the case of the discrete derivative matrix as in \eqref{dermtrx}, the computation of $D^\dagger_h$ can be done by means of rescaled fast Fourier transforms (indeed the matrix $D_h$ is almost a circulant matrix). On the other hand the discrete derivative matrix in 2D and in particular on meshes which are NOT on a regular square geometry (for example the one adapted to the discontinuity as in \cite{chdm99,boch00}) looses his structure and $D_h^\dagger$ is expected to be a fully populated matrix whose application and storage might be significantly demanding. How to deal with these situations is indeed a serious open problem. However, at the cost of a minor increase of implementation complexity, the use of subspace correction/domain decomposition strategy !
% ! as proposed in \cite{foDD} can turn out to offer a solution which is still computationally competitive. 
% \end{remark}

\subsubsection{The 2-D case, discrete Schwartz conditions, and constrained optimization}
\label{2dcase}
Let us assume now $d=2, \Omega = [0,1]^2$, and again $\alpha=\beta=1$.  Denote $u_{i,j}=u(h i,h j)$, $(i,j) \in \mathbb{Z}^2$, a discrete function defined on $\Omega_h:=\Omega \cap h \mathbb{Z}^2$, $n = \lfloor 1/h \rfloor$,  and
\begin{eqnarray*}
\mathcal J_{1/\sqrt{h}}(u ) &:=& h^2 \sum_{(h i, h j) \in \Omega_h} W_h \left ( \frac{u_{i+1,j} - u_{i,j}}{h} \right ) \\
&+& h^2 \sum_{ ( h i, h j) \in \Omega_h} W_h \left ( \frac{u_{i,j+1} - u_{i,j}}{h} \right )\\
&+& h^2 \sum_{(h i,h j) \in \Omega_h} ( u_{i,j} - g_{i,j})^2.
\end{eqnarray*} 
In two dimensions, we have to  consider the derivative matrix $D_h : \mathbb R^{n^2} \to \mathbb R^{2n(n-1)}$ that maps the vector $(u_{j + (i-1) n}) := (u_{i,j})$ to the vector composed of the finite differences in the horizontal and vertical directions $u_x$ and $u_y$ respectively, given by
$$
D_h u := \left [ 
\begin{array}{l}
u_x\\
u_y
\end{array}
\right ], \quad \left \{ \begin{array}{ll} (u_x)_{j + n(i-1) }:= (u_x)_{i,j}:= \frac{u_{i+1,j} - u_{i,j}}{h}, i=1,\dots,n-1, j=1,\dots,n\\
 (u_y)_{j + (n-1)(i-1) }:= (u_y)_{i,j}:= \frac{u_{i,j+1} - u_{i,j}}{h}, i=1,\dots,n, j=1,\dots,n-1
\end{array} \right. .
$$
{ Note that its range $R(D_h) \subset \mathbb{R}^{2n(n-1)}$ is a $(n^2-1)$-dimensional subspace} because $D_h c = 0$ for constant vectors $c \in \mathbb{R}^{n^2}$.  Again, we have the differentiation-integration formula, given by
$$
 u =D^\dagger_h D_h u + c,
$$
where $D_h^\dagger$ is the pseudo-inverse matrix of $D_h$ (in the Moore-Penrose sense); note that $D^\dagger_h$ maps $R(D_h)$ injectively into $\mathbb R^{n^2}$.  Also, $c$ is a constant vector that depends on $u$, and the values of its entries coincide with the mean value $h^2\sum_{ (hi,hj) \in \Omega_h} u_{i,j}$ of $u$.
\\
\\
Proceeding as before and again with a slight abuse of notation, we can reformulate the original discrete functional $\eqref{discrete}$ in terms of derivatives, and mean values, by
\begin{eqnarray*}
\mathcal J_{1/\sqrt h}(z,c) &=& h^2 \big[ \| D_h^\dagger z  - f \|^2_{\ell_2^{n^2}} + \| c - c_{g} \|_{\ell_2^{n^2}}^2 + \sum_{i,j} \min \left \{|z_{i,j}|^2, \frac{1}{h}  \right \} \big].
\label{2Dunconstrained}
\end{eqnarray*}
where $z = D_h u \in \mathbb{R}^{2n(n-1)}$, and $f =  D^{\dagger}_h D_h g \in \mathbb{R}^{n^2}$. Of course $c = c_g$ is again assumed at the minimizer $u$, since this latter term in $\mathcal J_{1/\sqrt h}$ does not depend on $z$.  However, in order to minimize only over vectors in $\mathbb R^{2n(n-1)}$ that are derivatives of vectors in $\mathbb R^{n^2}$, we must minimize $\mathcal J_{1/\sqrt{h}}(z,c)$ subject to { the constraint $D_h D^{\dagger}_h z = z$.}
\begin{figure}[htp]
\begin{center}
\includegraphics[width=2.5in]{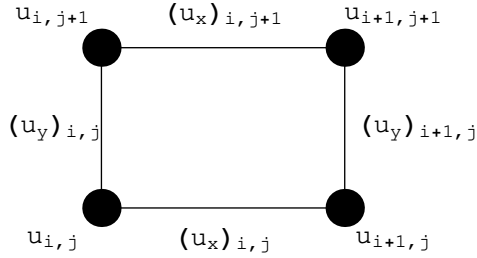}
\end{center}

\caption{Compatibility conditions of derivatives in 2D.}\label{compatib}
\end{figure}
\\
\\
The { $2n(n-1)$ linearly independent constraints $D_h D^{\dagger}_h z = z$} are equivalent to the \emph{discrete Schwartz constraints}\footnote{These discrete conditions correspond to the well-known  Schwartz mixed derivative theorem for which $\partial_{xy} u = \partial_{yx} u$ for any $u \in C^2(\Omega)$.}, 
\begin{equation}
\label{comp}
(u_y)_{i,j} + (u_x)_{i,j+1} =  (u_y)_{i+1,j} + (u_x)_{i,j},
\end{equation}
that establish the equivalence of the length of the paths from $u_{i,j}$ to $u_{i+1,j+1}$, whether one moves in vertical first and then in horizontal direction or in horizontal first and then in vertical direction (see Figure \ref{compatib}).
\\
\\
In short, we arrive at the following constrained optimization problem:

\begin{eqnarray}
&&
 \left\{ \begin{array}{llll}
\textrm{Minimize} &  \mathcal J_{1/\sqrt h}(z) = h^2 \big[ \| T z  - f \|^2_{\ell_2^{n^2}} +  \sum_{i,j} \min \left \{|z_{i,j}|^2, \frac{1}{h}  \right \} \big].
\\
\\
 \textrm{subject to} & {\cal Q} z = 0, \\\end{array}
\right.
\label{MS2d}
\end{eqnarray}
for $T = D_h^\dagger$ and { ${\cal Q} = {\cal I}- D_h D^{\dagger}_h $.}  Once the minimal derivative vector $z$ is computed, we can assemble the minimal $u$ by incorporating the mean value of $g$ as follows:
$$
 u = D_h^\dagger z + c_g.
$$

\subsubsection{Regularization of inverse problems by means of the Mumford-Shah constraint}
\label{MS4invprob}

The Mumford-Shah regularization term 
\begin{equation}
\label{MSreg}
MS ( u ) = \int_{\Omega \setminus S_u} | \nabla u |^2  +\beta \mathcal{H}^{d-1}(S_u),
\end{equation}
has been used frequently in inverse problems for image processing \cite{essh02,rari07}, such as inpainting and tomographic inversion. Despite the successful numerical results observed  in the aforementioned papers for the minimization of functionals of the type 
\begin{equation}
\label{regMS}
\mathcal J (u) = \alpha \| K u - g \|_{L^2(\Omega)} + MS(u),
\end{equation}
where $K:L^2 (\Omega) \to L^2 (\Omega)$ is a bounded operator which is not boundedly invertible, no rigorous results on existence of minimizers are currently available in the literature. Indeed, the Ambrosio compactness theorem \cite{AM89} used for the proof of the case $K = I$ does not apply in general. A few attempts towards using the regularization $MS$ for inverse problems in fracture detection appear in the work of Rondi \cite{ro07,ro08-2,ro08-1}, although restrictive technical assumptions on the admissible discontinuities of the solutions are required.
\\

\noindent As one of the contributions to this paper, we show that discretizations of regularized functionals of the type \eqref{regMS} \emph{always} have minimizers (see Theorem \ref{existmin}).  More precisely, these discretizations correspond to functionals of the form,
\begin{equation}
\mathcal J_{\sqrt{\beta/h}}(u ):=  \alpha h^2 \| K u - g\|_{\ell_2}^2 + h^2 \sum_{(h i, h j) \in \Omega_h} \left [ W_h \left ( \frac{u_{i+1,j} - u_{i,j}}{h} \right ) + W_h \left ( \frac{u_{i,j+1} - u_{i,j}}{h} \right ) \right ].
\label{regdiscrete}
\end{equation} 
and we prove that such functionals admit minimizers. Note that the discrete Mumford-Shah approximation $\eqref{discrete}$ can be written in this form.   We go on to show that such minimizers can be characterized by certain fixed point conditions, see Theorem 
 \ref{mintheorem} and Theorem \ref{globalmin}.  
As a consequence of these achievements we can prove that global minimizers are always isolated, although not necessarily unique, whereas local minimizers may constitute a continuum of unstable equilibria. Hence, our analysis will shed light on fundamental properties, virtues, and limitations, of regularization by means of the Mumford-Shah functional $MS$, and provide a rigorous justification of the numerical results appearing in the literature.
\\
\\
\noindent It is useful to show how the discrete functional \eqref{regdiscrete} can be still expressed in terms of the sole derivatives for general $K$. As done before in the case $K = {\cal I}$, and with the now usual identification $u= (D_h u, c)$, we can rewrite the functional in terms of derivatives and mean value as follows:
\begin{eqnarray}
\mathcal J_{\sqrt{\beta/h}}(z,c) &=&  h^2 \alpha \| K D_h^\dagger z  - (g - K c) \|_2^2+ h^2  \sum_{i,j} \min \left \{|z_{i,j}|^2, \frac{\beta}{h}  \right \},
\label{withc}
\end{eqnarray}
Note that in general we cannot anymore split orthogonally the discrepancy $\| K D_h^\dagger z  - (g - K c) \|_2^2$ into a sum of two terms which depend only on derivatives $z$ and mean value $c$ respectively. Nevertheless, for fixed $z$, it is straightforward to show that { $\bar{c} = \arg \min_c \mathcal J_{\sqrt{\beta/h}}(z,c)$  depends on $z$ via an affine map}. Indeed we can compute
$$
{\it \bar{c}} = \left ( \frac{ \langle K \mathbf 1, g - K D_h^\dagger z \rangle}{\| K \mathbf 1\|_{\ell^2}^2} \right ) \mathbf 1,
\label{c}
$$
where $\mathbf 1$ is the constant vector with entries identically $1$. Here we assume that $\mathbf 1 \notin \ker K$, that is a necessary condition in order to be able to identify the mean value of minimizers (a similar condition is required anytime we deal with regularization functionals which depend on the sole derivatives, see, e.g., \cite{chli97,ve01}).  By substituting this expression for $\bar{c}$ into \eqref{withc}, it is clear that the minimization of functionals $\eqref{regdiscrete}$ can be reformulated, in terms of the sole derivatives, as constrained minimization problems of the form $\eqref{MS2d}$.  

\section{Existence of minimizers for a class of discrete free-discontinuity problems}
In light of the observations above, we can transform the problem of the minimization of functionals of the type \eqref{regMS}, by means of discretization first and then reduction to sole derivatives, into the (possibly, but not necessarily) constrained minimization problem:
\begin{eqnarray}
&&
 \left\{ \begin{array}{llll}
\textrm{Minimize} &  \mathcal J_r (u) = \big[ \| T u  - g \|^2_{\ell_2^M} +  \sum^N_{i=1} \min \left \{|u_i|^2, r^2  \right \} \big].
\\
 \textrm{subject to} & {\cal Q} u = 0. \\\end{array}
\right.
\label{L2general}
\end{eqnarray}

Our first result ensures the existence of minimizers for the constrained optimization problem \eqref{L2general}:
\begin{prop}
Assume $r>0$, and fix linear operators $T:\mathbb R^N \to \mathbb R^M$ and ${\cal Q}: \mathbb R^N \to \mathbb R^{M'}$, which are identified in the following with their matrices with respect to the canonical bases. We also fix $g \in \mathbb R^M$. The constrained minimization problem
\begin{eqnarray}
&&
 \left\{ \begin{array}{llll}
\textrm{Minimize} &  \mathcal J_r (u) = \big[ \| T u  - g \|^2_{\ell_2^M} +  \sum^N_{i=1} \min \left \{|u_i|^2, r^2  \right \} \big]
\\
 \textrm{subject to} & {\cal Q} u = 0. \\\end{array}
\right.
\label{L2}
\end{eqnarray}
has minimizers $u^*$. 
\label{existmin2}
\end{prop}
\begin{proof}
We begin by noting that $\inf_{{\cal Q}u = 0} {\cal{J}}_r({u})$ is well-defined and finite, since ${\cal{J}}_r \geq 0$ is bounded from below.  It remains to show that there exists a vector ${u^*}$ that satisfies ${\cal{J}}_r({u^*}) = \inf_{u \in \mathbb R^N} {\cal{J}}_r({u})$.  Towards this goal, consider the following partition $\mathcal P = \{ \mathcal U_{\mathcal I_j } \}_{j=1}^{2^N}$ of $\mathbb{R}^N$ indexed by the subsets $\mathcal I_j$ of the index set $\mathcal I = \{1,2,...,N\}$, as follows:
\begin{equation}
\mathcal U_{{\cal I}_j} := \{u \in \mathbb{R}^{N}: |u_i| \leq r, i \in {\cal I}_j, |u_i| > r, i \in {\cal I} / {\cal I}_j \}.
\end{equation}
The minimization  of ${\cal{J}}_r$ subject to ${\cal Q}u = 0$ and constrained to the closure of the subset $\mathcal U_{\mathcal I_j}$ can be reformulated as a quadratic optimization problem, for which the classical Frank-Wolfe theorem  $\cite{bashsh93}$ guarantees the existence of a minimizer $u(\mathcal I_j)$.  Now, since $\mathbb{R}^N = \cup_j {\cal I}_j$, the minimal value of $\mathcal J_r$ subject to $Qu = 0$ and over all of $\mathbb{R}^N$ is just the minimal value from the finite set $\{ \mathcal J_r(u({\cal I}_j)): \quad j=1, \dots, 2^N\}$; that is,
\begin{eqnarray}
\min_{{\cal Q}u = 0} {\cal{J}}_r({u}) &=& \min_{{\cal I}_j \subset {\cal I}} {\cal{J}}_r({u}({\cal I}_j))
\nonumber
\end{eqnarray}
and $u^*=\arg \min_{{\cal Q}u = 0} {\cal{J}}_r({u}) = u\big(\arg \min_{{\cal I}_j \subset {\cal I}} {\cal{J}}_r({u}({\cal I}_j))\big).$
\end{proof}
In fact, Proposition \ref{existmin2} extends to a much larger class of free-discontinuity type minimization problems; by the same reasoning as before, we arrive at the more general result:
\begin{theorem}
The constrained minimization problem
\begin{eqnarray}
&&
 \left\{ \begin{array}{llll}
\textrm{Minimize} &  \mathcal J^p_r (u) = \big[ \| T u  - g \|^2_{\ell_2^M} +  \sum^N_{i=1} \min \left \{|u_i|^p, r^p  \right \} \big]
\\
 \textrm{subject to} & {\cal Q} u = 0. \\\end{array}
\right.
\label{Lpgeneral}
\end{eqnarray}
has minimizers $u^*$ for any real-valued parameter $p \geq 1$. 
\label{existmin}
\end{theorem}
The Frank-Wolfe theorem, which guarantees the existence of minimizers for quadratic programs with bounded objective function, does not apply to the general case $p \geq 1$ where the objective function $\mathcal J^p_r$ is not necessarily quadratic.  Nevertheless, with the following generalization for the Frank-Wolfe theorem, Theorem \ref{existmin} follows directly from a similar argument as for Proposition \ref{existmin2}.

\begin{prop}
Suppose $A$ is an $N \times N$ positive semidefinite matrix, and suppose $b$ and $c$ are $N \times 1$ vectors.   Suppose also that $X$ is a nonempty convex polyhedral subset of $\mathbb{R}^N$.  The convex optimization problem
\begin{eqnarray}
 \left\{ \begin{array}{llll}
\textrm{minimize} & u^t A u + b^t u + \sum_{1 \leq j \leq N} c_j |u_j|^p
\\
 \textrm{subject to} & u \in X.  \\\end{array} 
\right.
\label{FWgen}
\end{eqnarray}
admits minimizers for any real parameter $p \geq 1$, as long as the objective function is bounded from below.  
\label{firstlemma}
\end{prop}
For ease of presentation, we reserve the proof of Proposition \ref{firstlemma} to the Appendix.

\noindent From the proof of Theorem \ref{existmin}, one could in principle obtain a minimizer for ${\cal J}^p_r$ by computing a minimizer ${u}({\cal I}_j)$ for each subset ${\cal I}_j \subset {\cal I}$ using a quadratic program solver $\cite{bashsh93}$, and then minimizing ${\cal J}^p_r$ over the finite set of points $\{ u({\cal I}_j) \}$.  Unfortunately, this algorithm is computationally infeasible as the number of subsets of the index set $\{1,2, ..., N \}$ grows exponentially with the dimension $N$ of the underlying space.   
Indeed, the minimization problem \eqref{Lpgeneral} is  NP-hard, as the known NP-complete problem SUBSET-SUM can be reduced to this problem. A complete discussion about the NP-hardness of \eqref{Lpgeneral} can be found in \cite{RachelPrep}. 

%{\mnew

\section{An iterative thresholding algorithm for 1-D free-discontinuity inverse problems}

\subsection{Overview of the algorithm}

In this section, we introduce an algorithm that is guaranteed to converge to a local minimizer of the real-valued functional ${\cal{J}}_r^p: \ell_2(\mathcal I) \rightarrow \mathbb{R}$ having the form 
\begin{equation}
{\cal{J}}^p_r ({ u}) =  \| T {u} - g \|^2_{\ell_2(\mathcal K)} +  \sum_{i \in \mathcal I} \min \{ |u_i|^p, r^p \},
\label{generalform}
\end{equation}
subject to the { conditions:
\begin{itemize}
\item  $\mathcal I$ and $\mathcal K$ are countable sets of indices, and $T:\ell_2(\mathcal I) \to  \ell_2(\mathcal K)$ is a bounded linear operator, which is in the following identified with its matrix associated to the canonical basis;
\item the operator $T$ has spectral norm $\|T\|  < 1$. Note that this requirement is easily met by an appropriate scaling for the functional, i.e., we may have to consider instead
$$
{\cal{J}}^p_r ({ u}) =  \gamma \| T{ u} - g \|^2_{\ell_2(\mathcal K)} +  \gamma \sum_{i \in \mathcal I} \min \{ |u_i|^p, r^p \}, \quad \gamma \leq 1.
$$ 
This modification leads to minor changes in the analysis that follows (see also Subsection \ref{denoising}), and throughout this paper we assume, without loss of generality, that $\gamma =1$;

\item  the parameter $p$ is in the range $1 \leq p \leq 2$.  In case the index set $\mathcal I$ is \emph{finite}, only the restriction $p \geq 1$ is necessary.  

\end{itemize}
We note that the scaled 1D discrete Mumford-Shah functional $\frac{1}{h}{\cal{J}}_{1/\sqrt h}$ is clearly a functional of the form \eqref{generalform} having $r=1/\sqrt h$, index set $\mathcal I  = \{1, \dots, \lfloor r^2 \rfloor\}$, parameter $p = 2$, and operator $T = D_{1/r^2}^{\dagger}: \mathbb{R}^{\lfloor r^2 \rfloor - 1} \rightarrow \mathbb{R}^{\lfloor r^2 \rfloor}$ .  As shown in the Appendix, the operators $D_{1/r^2}^{\dagger}$ satisfy the uniform bound $\| D_{1/r^2}^{\dagger} \| \leq 1/2$, independent of dimension, so a scaling factor is not needed in this case.  
\\

\noindent In the following, we will not minimize ${\cal{J}}^p_r$ directly.  Instead, we propose a \emph{majorization-minimization} algorithm for finding solutions to ${\cal{J}}^p_r$, motivated by the recent application of such algorithms for minimizing energy functionals arising in sparse signal recovery and image denoising \cite{blda??,DDD}.    More precisely, consider the following \emph{surrogate} objective function,
\begin{equation}
 {\cal{J}}^{p, surr}_r (u,a) := {\cal{J}}^p_r(u) - \| T{u} - T {a} ||_{\ell_2(\mathcal K)}^2 + \| {u} - {a} \|_{\ell_2(\mathcal I)}^2. \quad u,a \in \ell_2(\mathcal I).
\label{surr1}
 \end{equation}
The surrogate functional ${\cal{J}}^{p,surr}_r$ satisfies ${{\cal{J}}}^{p,surr}_r (u, {a}) \geq {\cal{J}}^p_r({ u})$ everywhere, with equality if and only if ${u = a}$, and is such that the sequence
\begin{equation}
{ u}^{n+1} = \arg \min_{ u} {{\cal{J}}}^{p,surr}_r ({ u}, { u}^{n})
\label{map}
\end{equation}
obtained by successive minimizations of ${{\cal{J}}}^{p,surr}_r ({ u, a})$ in ${ u}$ for fixed ${ a}$ results in a nonincreasing sequence of the original functional ${\cal{J}}^p_r (u^n)$ (see Lemmas \ref{l1} and \ref{l2}).   We will study the implementation and the convergence properties of the iteration \eqref{map} as follows:
\begin{itemize}
\item in Section $3.2$, we review the standard properties of majorization-minimization iterations,
\item in Section $3.3$, we explicitly compute $u$-global minimizers of the surrogate functional ${{\cal{J}}}^{p,surr}_r ({ u, a})$, for $a$ fixed; 
\item in Section $3.4$ we discuss a connection between the resulting thresholding functions and thresholding functions used in sparse recovery,
\item in Sections $3.5$, $3.6$, and $3.7$, we show that  the sequence $({ u}^n )_{n \in \mathbb N}$ defined by \eqref{map} will converge to a stationary value $\bar{u}= \arg \min_{ u} {{\cal{J}}}^{p,surr}_r ({ u}, \bar{ u})$, starting from any initial value ${ u}^0$ for which ${\cal{J}}^p_r({ u}^0) < \infty$,
\item in Section $3.8$, we show that such stationary values $\bar{u}$ are also local minimizers of the original functional ${\cal{J}}^p_r$ that satisfy a certain fixed point condition, and 
\item in Section $3.9$, it is shown that any global minimizer of ${\cal{J}}^p_r$ is among the set of possible fixed points $\bar{u}$ of the iteration $\eqref{map}$.
\end{itemize}
By means of the thresholding algorithm, we also show that global minimizers of the functional ${\cal{J}}^p_r$ are isolated, and moreover possess a certain segmentation property that is also shared by fixed points of the algorithm.
}

\subsection{Preliminary lemmas}
The lemmas in this section are standard when using surrogate functionals (see $\cite{DDD}$ and $\cite{blda??}$), and concern general real-valued surrogate functionals of the form
\begin{equation}
{\cal{F}}^{surr}({ u, a}) = {\cal{F}}({ u}) - \|T{ u} - T{ a}\|^2_{\ell_2(\mathcal K)}+ \|{ u - a}\|^2_{\ell_2(\mathcal I)}.
\label{surrgen}
\end{equation}
The lemmas in this section hold independent of the specific form of the functional ${\cal{F}}:\ell_2(\mathcal I) \to \mathbb R^+$,  but do rely on the restriction that $\|T\| < 1$.  
\begin{lemma}
If the real-valued functionals ${\cal{F}}({ u})$ and ${\cal{F}}^{surr}({ u,a})$ satisfy the relation $\eqref{surrgen}$ and the sequence $({ u^n})_{n \in \mathbb N}$ defined by ${ u}^{n+1} = \arg \min_{ u \in \ell_2(\mathcal I)} {\cal{F}}^{surr} ({ u, u^n})$ is initialized in such a way that ${\cal{F}}({ u^0}) < \infty$, then the sequences ${\cal{F}} ({ u^n})$ and ${{\cal{F}}}^{surr} ({ u^{n+1}, u^{n}})$ are non-increasing as long as $\|T\| < 1$.
\label{l1}
\end{lemma}
\begin{proof}
Since $\|T\| < 1$, also $\|T^*T\| < 1$,  and so the operator $L = \sqrt{I - T^*T}$ is a well-defined positive operator whose spectrum is contained within a closed interval $[c, 1]$ that is bounded away from zero $c > 0$.   We can then rewrite ${{\cal{F}}}^{surr} ({ u^{n+1} u^{n}})$ as ${{\cal{F}}}^{surr} ({ u^{n+1}, u^{n}}) = {{\cal{F}}}({ u^{n+1}}) + \|L({ u^{n+1} - u^n)}\|_{\ell_2(\mathcal I)}^2$, from which it follows that
\begin{eqnarray}
{{\cal{F}}}({ u^{n+1}}) &\leq& {{\cal{F}}}({ u^{n+1}}) +  \|L({ u^{n+1} - u^n})\|_{\ell_2(\mathcal I)}^2 \nonumber \\
&=& {{\cal{F}}}^{surr} ({ u^{n+1}, u^{n}}) \nonumber \\
&\leq& {{\cal{F}}}^{surr} ({ u^{n}, u^{n}}) \nonumber \\
&=& {{\cal{F}}}({ u^{n}}) \nonumber \\
&\leq& {{\cal{F}}}({ u^{n}}) +  \|L({ u^{n} - u^{n-1}})\|_{\ell_2(\mathcal I)}^2 \nonumber \\
&=& {{\cal{F}}}^{surr} ({ u^{n}, u^{n-1}}),
\label{est}
\end{eqnarray}
where the second inequality follows from ${ u^{n+1}}$ being a minimizer of $ {{\cal{F}}}^{surr} ({ u, u^{n}})$.
\end{proof}
\noindent From Lemma \ref{l1} we obtain the following corollary:
\begin{lemma}
As long as the conditions of Lemma $\ref{l1}$ are satisfied, one can choose $N \in \mathbb N$ sufficiently large such that for all $n \geq N$, $\|{ u^{n+1} - u^n} \|_{\ell_2(\mathcal I)} \leq \epsilon$, i.e., $$\lim_{n \to \infty}  \|{ u^{n+1} - u^n} \|_{\ell_2(\mathcal I)}=0.$$
\label{l2}
\end{lemma}
\begin{proof}
From Lemma  \ref{l1}, it follows that ${\cal F} (u^{n}) \geq 0$ is a nonincreasing sequence, therefore it converges, and ${\cal F} (u^{n}) - {\cal F} (u^{n+1}) \to 0$ for $n \to \infty$. The lemma follows from \eqref{est}, and the estimates
$$
{\cal F} (u^{n}) - {\cal F} (u^{n+1}) \geq \|L({ u^{n+1} - u^n})\|_{\ell_2(\mathcal I)}^2 \geq (1 - \|T\|^2) \|  u^{n+1} - u^n\|_{\ell_2(\mathcal I)}^2.
$$
\end{proof}

% Lemma $\eqref{l2}$ follows if we can show that the sequence $s_N = \sum_{n=0}^N ||{ u^{n+1} - u^n} ||_2^2$ converges as $N \rightarrow \infty$.  Since the terms $||{ u^{n+1} - u^n} ||_2^2$ are nonnegative, the sequence $( s_N ) $ is monotonically increasing, and thus converges if and only if the $s_N$ are uniformly bounded.  But        

% \begin{eqnarray}
% \sum_{n=0}^N ||{ u^{n+1} - u^n} ||_2^2 &=& \sum_{n=0}^N || L^{-1}L({ u^{n+1} - u^n)} ||_2^2 \nonumber \\
% &\leq& \frac{1}{c} \sum_{n=0}^N || L({ u^{n+1} - u^n}) ||_2^2 \nonumber \\
% &\leq& \frac{1}{c} \sum_{n=0}^N \big[ {{\cal{F}}} ({ u^{n}}) - {{\cal{F}}} ({ u^{n+1}}) \big] \nonumber \\
% &=& \frac{1}{c} ( {{\cal{F}}} ({ u_{0}}) - {{\cal{F}}} ({ u^{n+1}}) ) \nonumber \\
% &\leq& \frac{1}{c} {{\cal{F}}} ({ u_{0}});
% \end{eqnarray}
% Again, the constant $c$ is a positive lower bound on the spectrum of the positive operator $L$.  The second inequality uses Lemma $\eqref{l1}$.

\subsection{The surrogate functional ${{\cal{J}}}^{p, surr}_r$, its explicit minimization, and a new thresholding operator}
It is not immediately clear that the surrogate functional ${{\cal{J}}}^{p,surr}_r$ in $\eqref{surr1}$ is any easier to manage than its parent functional ${{\cal{J}}}^p_r$.  However, expanding the squared terms on the right hand side of $\eqref{surr1}$, ${{\cal{J}}}^{p,surr}_r(u,a)$ can be equivalently expressed as
\begin{eqnarray}
{{\cal{J}}}^{p,surr}_r ({ u, a}) &=& \| { u} - (I - T^*T){ a} + T^*{ g}  \|_{\ell_2(\mathcal I)}^2 + \sum_{i \in \mathcal I} \min \{ |u_i|^p, r^p \} + C\nonumber \\ 
&=& \sum_{i \in \mathcal I} \Big[ ( u_i - [a - T^*Ta + T^*g]_i )^2 + \min \{ |u_i|^p, r^p \} \Big] + C, \nonumber
\label{surr}
\end{eqnarray}
where the term $C = C(T,a,g)$ depends only on $T$, ${ a}$ and ${ g}$.  Indeed, unlike the original functional ${\cal{J}}^p_r$, the surrogate functional ${{\cal{J}}}^{p,surr}_r$ {\it decouples} in the variables $u_i$, due to the cancellation of terms involving $\|T{ u}\|_{\ell_2}^2$.  Because of this decoupling, global ${ u}$-minimizers of ${{\cal{J}}}^{p,surr}_r({ u,a})$, for $a$ fixed, can be computed {\it component-wise} according to
\begin{equation}
\bar{u}_i = \arg \min_{t \in \mathbb R}  \Big[ (t - [a - T^*Ta + T^*g]_i )^2 + \min \{ |t|^p, r^p \} \Big], \quad i \in \mathcal I.
\label{surmin}
\end{equation} 
One can solve $\eqref{surmin}$ explicitly when e.g. $p = 2$, $p = 3/2$, and $p = 1$; in the general case $p \geq 1$, we have the following result: 
\begin{prop}[Minimizers of ${{\cal{J}}}^{p,surr}_r({ u,a})$ for $a$ fixed]\label{thmthrs}.
\begin{enumerate}
\item If {$p > 1$}, the minimization problem ${ \bar{u}} = \arg \min_{ u \in \ell_2(\mathcal I)} {{\cal{J}}}^{p,surr}_r ({ u}, { a})$ can be solved component-wise by
\begin{equation}
\bar{u}_i = H_{(p,r)}([{ a} - T^*T{ a} + T^*{ g}]_i), \quad i \in \mathcal I,
\label{three}
\end{equation}
where $H_{(p,r)}: \mathbb{R} \rightarrow \mathbb{R}$ is the \emph{`thresholding function'},
\begin{equation}
H_{(p,r)} (\lambda) =
\left\{
\begin{array}{ll}
F_p^{-1}(\lambda), & |\lambda| \leq \lambda'(r,p)  \\
\lambda, & |\lambda|  > \lambda'(r,p). \\  
\end{array}
\right.
\label{thresh}
\end{equation}
Here, $F_p^{-1}(\lambda)$ is the inverse of the function $F_p(t) = t + \frac{p}{2}\sgn{t}{|t|}^{p-1}$,
and $\lambda':=\lambda'(r,p) \in (r, r + \frac{p}{2}r^{p-1})$ is the unique positive value at which
\begin{equation}
(F_p^{-1}(\lambda') - \lambda')^2 + |F_p^{-1}(\lambda')|^p = r^p.
\end{equation}
\item When {  $p=1$}, the general form $\eqref{three}$ still holds, but we have to consider two cases:
\begin{enumerate}
\item If $r > 1/4$, the thresholding function $H_{(1,r)}: \mathbb{R} \rightarrow \mathbb{R}$ satisfies
\begin{equation}
H_{(1,r)} (\lambda) =
\left\{
\begin{array}{ll}
0, & |\lambda| \leq 1/2  \\
(|\lambda| - 1/2)\sgn{\lambda}, &1/2 < |\lambda|  \leq r + 1/4 =  \lambda'(r,1) \\
\lambda, & |\lambda| >  r + 1/4
\end{array}
\right.
\label{thresh}
\end{equation}
\item If, on the other hand, $r \leq1/4$, the function $H_{(1,r)}$ satisfies
\begin{equation}
H_{(1,r)} (\lambda) =
\left\{
\begin{array}{ll}
0, & |\lambda| \leq \sqrt{r}=  \lambda'(r,1)  \\
\lambda, & |\lambda| > \sqrt{r}
\end{array}
\right.
\label{thresh}
\end{equation}
\end{enumerate}
\end{enumerate}
In all cases, the function $H_{(p,r)}$ is continuous except at $\lambda'(r,p)$, where $H_{(p,r)}$ has a jump-discontinuity of size $\delta(r,p) = |\lambda' - H_{(p,r)}(\lambda')| > 0$ if $r > 0$.  In particular, it holds that $\lambda'(r,p) > r$ while $H_{(p,r)}(\lambda') < r$.
\end{prop}
\noindent We leave the proof of Proposition \ref{thmthrs} to the Appendix.  

\begin{figure}[htp]
\begin{center}
\subfigure[]{ \label{1(a)}
\includegraphics[width=3in]{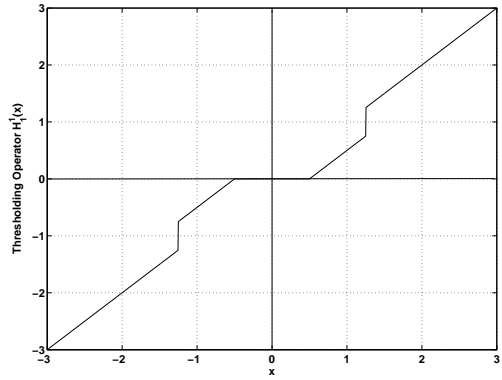}}
\subfigure[]{\label{1(b)}
\includegraphics[width=3in]{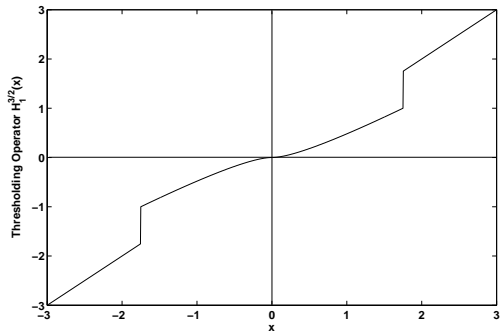}}
\subfigure[]{\label{1(c)}
\includegraphics[width=3in]{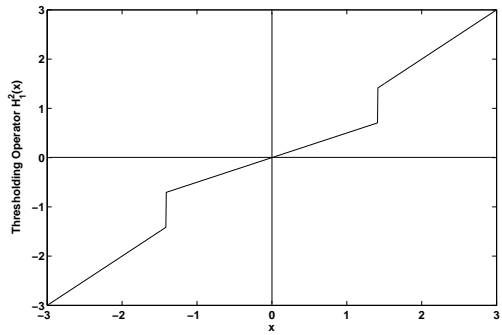}}
\end{center}
\caption{The discontinuous thresholding functions  $H_{(1,1)}$,  $H_{(3/2,1)}$, and $H_{(2,1)}$, with parameters $p=1,3/2$, and $2$, respectively, and $r = 1$.}\label{threshold}
\end{figure}

\begin{remark}
\emph{ In the particular case $p = 2$ corresponding to classical Mumford-Shah regularization \eqref{L2general}, the thresholding function $H_{(2,r)}: \mathbb{R} \rightarrow \mathbb{R}$ has a particularly simple explicit form:}
\begin{equation}
H_{(2,r)} (\lambda) =
\left\{
\begin{array}{ll}
\lambda/2, & |\lambda| \leq \sqrt{2}r  \\
\lambda, & |\lambda|  > \sqrt{2}r \\  
\end{array}
\right. 
\label{thresh}
\end{equation}
\emph{In addition to $H_{(2,r)}$ and $H_{(1,r)}$, the thresholding operator $H_{(3/2,r)}(\lambda)$ corresponding to $p = 3/2$ can also be computed explicitly, by solving for the positive root of a suitable polynomial of third degree.  In Figure $\ref{threshold}$ below, we plot $H_{(2,1)}, H_{(3/2,1)}$, and $H_{(1,1)}$ with parameter $r = 1$.  For general noninteger values of $p$, $H_{(p,r)}$ cannot be solved in closed form.   However, recall the following general properties of $H_{(p,r)}$:
\begin{itemize}
\item $H_{(p,r)}$ is an odd function, 
\item $H_{(p,r)}(0) = 0$, and 
\item $H_{(p,r)}(\lambda) = \lambda$ once $|\lambda| > r + \frac{p}{2}r^{p-1}$. 
\end{itemize}
In fact, we can effectively \emph{pre}compute $H_{(p,r)}$ by numerically solving for the value of $H_{(p,r)}(\lambda_j)$ on a discrete set $\{\lambda_j\}$ of points $\lambda_j  \in (0,\frac{p}{2}r^{p-1} + r]$.  At $\lambda_j$, one just needs to solve the real equation
\begin{equation}
h_j + \frac{p}{2}h_j^{p-1} - \lambda_j = 0
 \end{equation}
which can be computed effortlessly via a root-finding procedure such as Newton's method:  while $h_j$ satisfies $(h_j - \lambda_j)^2 + (h_j)^p \leq r^p$, set $H_{(p,r)}(\lambda_j) = h_j$; once this constraint is violated, set $H_{(p,r)}(\lambda_j) = \lambda_j$. }
\end{remark}

\subsection{Connection to sparse recovery}

When $p=1$ and $r \leq 1/4$, we know from Theorem \ref{thmthrs} that the iterative algorithm
\begin{equation}
{ u}^{n+1} = \arg \min_{ u} {{\cal{J}}}^{p,surr}_r ({ u}, { u}^{n})
\end{equation}
reduces to the component-wise thresholding
{\nnew 
\begin{eqnarray}\label{hard}
u^{n+1}_i &=& H_{\sqrt{r}}([u^{n}- T^*T{ u^{n}} + T^*g]_i),
\end{eqnarray}
}
where 
\begin{equation}
H_{\gamma}(\lambda) =
\left\{
\begin{array}{ll}
0, & |\lambda| \leq \gamma  \\
\lambda, & |\lambda|  > \gamma.
\end{array} \right.
\end{equation}
This thresholding function $H_{\gamma}: \mathbb{R} \rightarrow \mathbb{R}$ is referred to as {\it hard-thresholding} in the area of sparse recovery, and the iteration \eqref{hard} generated by successive applications of hard thresholding has been previously studied \cite{blda??}.  In particular, the iteration \eqref{hard} was shown in \eqref{hard} to correspond to successive minimization in $u$ for fixed $a$ of the surrogate functional ${\cal F}_r^{0,surr} (u, a)$ corresponding to the $\ell_0$ regularized functional, 
\begin{equation}\label{ell_0}
 \mathcal F^0_r(u) = \|T{ u} - { g}\|^2_{\ell_2(\mathcal K)} + r \|{ u}\|_{\ell_0(\mathcal I)}.
\end{equation}
Here, the $\ell_0$ quasi-norm  $\| u \|_{\ell_0(\mathcal I)} := \sum_{i \in \mathcal I } | u_i|_0$  is defined component-wise by $$
|u_i|_0= \left \{ 
\begin{array}{ll}
0, & \textrm{ if } u_i= 0\\
1, & \textrm{ otherwise} 
\end{array} \right . 
$$  
The $\ell_0$ regularized functional $\mathcal F^0_r(u)$ is related to the so-called \emph{K-sparse} problem,
\begin{eqnarray}\label{Ksparse}
&& \left\{ \begin{array}{llll}
\textrm{minimize} &  \| T u  - g \|_{\ell_2}^2 \\
 \textrm{subject to} & \|u \|_0 \leq K, 
\end{array}
\right.
\end{eqnarray}
in that there exists a $r$, that depends on $g$ and $K$, such that the solution to the $K$-sparse problem is the minimizer of the $\ell_0$ regularized functional. 
The $K$-sparse problem \eqref{Ksparse} is NP-hard in general \cite{RachelPrep}, but under certain restrictions on the matrix $T$, it is possible to solve \eqref{Ksparse} using fast algorithms.   For example, if the $m \times N$ matrix $T$ satisfies a certain \emph{restricted isometry property} of order $2K$ \cite{Candes08}, and there exists a $K$-sparse vector satisfying the constraint $T u = g$, then $u$ is the unique solution to \eqref{Ksparse} and can be recovered as the limit of the following {\nnew {\it iterative hard thresholding} (IHT) \cite{IHT}:
\begin{eqnarray}\label{hard_sparse}
u^{n+1} &=& \tilde {\mathbb{H}}_{K}(u^n - T^*T{ u^n} + T^*g).
\end{eqnarray}
Here, the thresholding operator $\tilde {\mathbb H}_s(u)$} sets all but the largest (in magnitude) $s$ elements of $u$ to zero.   This algorithm can be viewed as a variant of the hard thresholding algorithm \eqref{hard} with threshold parameter $r = r_n$ adaptively adjusted at each iteration to remain consistent with the knowledge that a $K$-sparse solution exists.   In fact, a modified version of IHT, called {\nnew \emph{normalized iterative hard thresholding} (NIHT), represents the state of the art among a large class of algorithms that have been designed to solve the $K$-sparse problem \eqref{Ksparse} under RIP or related assumptions on the matrix $T$ \cite{NIHT}, see also the paper repository \cite{CS.com}.}
Preliminary numerical results indicate that the performance of NIHT could be strengthened by replacing hard thresholding with a hybrid soft-hard thresholding, as shown at the top of Figure \ref{threshold}, as derived in Proposition \ref{thmthrs} from the minimization of free-discontinuity functional ${\cal J}_r^p$ with parameters $p = 1$ and $r > 1/4$.  
%Minimizers of the $K$-sparse problem are desired in the areas of sparse recovery and compressed sensing, when one is interested in solving inverse problems with sparse vector representations $\cite{carota06}$.  It was with this motivation that convergence of the sequence $\eqref{hard}$ to a local minimizer of the $\ell_0$ regularized functional $\eqref{ell_0}$ was shown in \cite{blda??}.   According to Proposition \ref{thmthrs}, and Theorem \ref{mintheorem} below, this \emph{very same} thresholding algorithm also converges to local minimizers of the functional, 
%\begin{equation}
%{\cal{J}}^1_r ({ u}) =  \| T{ u} - { g} \|^2_{\ell_2(\mathcal K)} +  \sum_{i \in \mathcal I} \min \{ |u_i|, r \},
%\label{ell1-free}
%\end{equation}
%suggesting an equivalence between minimizers of the free-discontinuity type functional $\eqref{ell1-free}$ and the sparsity-regularized functional $\eqref{ell_0}$.   
\\

Because a convergence analysis of the iteration $\eqref{hard}$ corresponding to hard thresholding  has been studied already $\cite{blda??}$, we omit the case $p = 1$ and $r \leq 1/4$ in the sequel.

%, but return its inclusion in section 3 when we discuss the connection between stationary points of the map $\eqref{map3}$ and local minimizers of the functionals $\eqref{generalform}$. 

\subsection{Fixation of the discontinuity set}
We prove now that the sequence $({ u}^n)_{n \in \mathbb N}$ defined by
\begin{eqnarray}
{ u}^{n+1} &=& \arg \min_{ u} {{\cal{J}}}^{p,surr}_r ({ u}, { u}^{n}) 
\end{eqnarray} 
or equivalently, according to Proposition $\ref{thmthrs}$, component-wise by
\begin{eqnarray}\label{map3}
u^{n+1}_i &=& H_{(p,r)}([u^n - T^*Tu^n + T^*g]_i), \quad i \in \mathcal I,
\end{eqnarray}
will converge, granted that $p \geq 1$ and $\|T\| < 1$.  To ease notation, we define the operator $\mathbb{H}: \ell_2(\mathcal I) \rightarrow \ell_2(\mathcal I)$ by its component-wise action,
\begin{equation}
[\mathbb{H}(u)]_i := H_{(p,r)}([u - T^*T{ u} + T^*g]_i);
 \label{u}
\end{equation}
so that the iteration $\eqref{map3}$ can be written more concisely in operator notation as
\begin{equation}
u^{n+1} = \mathbb{H}(u^n).
\end{equation}
We omit the dependence of $\mathbb{H}$ on the parameters $p, r$, and the function ${ g}$ for continuity of presentation.  At the core of the convergence proof is the fact that the `discontinuity set', indicated below by ${\cal I}_1^n$, of $u^n$ must eventually fix during the iteration \eqref{map3}, at which point the `free-discontinuity' problem is transformed into a simpler `fixed-discontinuity' problem. 

\begin{lemma}[Fixation of the index set ${\cal I}_1$]
Fix $p \geq 1, r \in \mathbb{R}^+$, and ${ g} \in \ell_2(\mathcal K)$.  Consider the iteration
\begin{equation}
u^{n+1} = \mathbb{H}(u^n)
\end{equation}
and the time-dependent partition of the index set $\mathcal I$ into `small' set
\begin{eqnarray}
\mathcal I_0^n &=&  \{ i \in \mathcal I: |u^n_i| \leq \lambda'(r,p) \} 
\label{lam0}
\end{eqnarray}
and 'large' set
\begin{eqnarray}
\mathcal I_1^n &=& \{ i  \in \mathcal I:  |u^n_i| > \lambda'(r,p) \} 
\label{lam1}
\end{eqnarray}
where $\lambda'(r,p)$ is the position of the jump discontinuity of the thresholding function, as defined in Proposition $\ref{thmthrs}$.  For $N \in \mathbb N$ sufficiently large, this partition fixes during the iteration $u^{n+1} = \mathbb{H} (u^n)$; that is, there exists a set $\mathcal I_0$ such that for all $n \geq N$, $\mathcal I_0^n = \mathcal I_0$ and $\mathcal I_1^n = \mathcal I_1 := \mathcal I \setminus \mathcal I_0$.  
\label{l3}
\end{lemma}

\begin{proof}
By {\it discontinuity} of the thresholding operator $H_{(p,r)}(\lambda)$, each sequence component
\begin{equation}
u^{n}_i = H_{(p,r)} ([u^{n-1} - T^*Tu^{n-1} + T^*g]_i)
\end{equation}
satisfies 
\begin{enumerate}
\item[(a)] $| u^n_i | \leq \lambda'(r,p) - \delta(r,p) < \lambda'(r,p) $, if $i \in \mathcal I_0^n$, or
\item[(b)] $| u^n_i | > \lambda'(r,p)$, if $i \in \mathcal I_1^n$.
\end{enumerate}
Thus, $|u^{n+1}_i - u^n_i| \geq \delta(r,p)$ if $i \in \mathcal I_0^{n+1} \bigcap \mathcal I_1^n$, or vice versa if $i \in \mathcal I_0^{n} \bigcap \mathcal I_1^{n+1}$.  At the same time, Lemma $\ref{l2}$ implies
\begin{equation}
| u^{n+1}_i - u^n_i | \leq \|u^{n+1} - u^n \|_{\ell_2(\mathcal I)} \leq \epsilon, 
\label{contra}
\end{equation}   
once $n \geq N(\epsilon)$, and $\epsilon > 0$ can be taken arbitrarily small.  In particular, $\eqref{contra}$ implies that $\mathcal I_0$ and $\mathcal I_1$ must be fixed once $n \geq N(\epsilon)$ and $\epsilon < \delta(r,p)$.
\end{proof}

\noindent After fixation of the index set $\mathcal I_0 = \{i \in \mathcal I: |u_i^n| \leq \lambda'(r,p)\}$, $\mathbb{H}(u^n) =\mathbb{U}_{\mathcal I_0}(u^n)$ and $\mathbb{U}_{\mathcal I_0}$ is an operator having component-wise action{, for $p>1$,}
\begin{eqnarray}
[\mathbb{U}_{\mathcal I_0} u]_i &=& \left\{
\begin{array}{ll} F_p^{-1}([(I - T^*T)u + T^*g]_i), & \textrm{if } i \in \mathcal I_0 \\  
\big( (I - T^*T)u + T^*g \big)_i, & \textrm{if } i \in \mathcal I_1
\end{array} \right.
\label{separate}
\end{eqnarray}

Here, as in Proposition $\ref{thmthrs}$, the function $F_p^{-1}$ is the inverse of the function $F_p(t) = t + \frac{p}{2}\sgn{t}|t|^{p-1}$.  Again, for ease of presentation, we omit the dependence of $\mathbb{U}_{\mathcal I_0}$ on the parameters $p,r$, and ${ g}$. { For $p=1$ the description is similar, and in general,} one easily verifies the equivalence
\begin{eqnarray}
\mathbb{U}_{\mathcal I_0} ({ v}) =  \arg \min_{ u \in \ell_2(\mathcal I)} {\cal{J}}_{\mathcal I_0}^{p,surr} ({ u, v})
\label{equiv}
\end{eqnarray}
where ${\cal{J}}_{\mathcal I_0}^{p,surr}$ is a surrogate for the {\it convex} functional,
\begin{equation}
 {\cal{J}}_{\mathcal I_0}^p({ u}) :=  \| T{ u} - { g} \|_{\ell_2(\mathcal K)}^2 +  \sum_{i \in \mathcal I_0} |u_i|^p. 
\end{equation}
That is, fixation of the index set $\mathcal I_0$ implies that the sequence $( u^n)_{n \in \mathbb N}$ has become constrained to a subset of $\ell_2(\mathcal I)$ on which the map $\mathbb{H}$ agrees with a map $\mathbb{U}_{\mathcal I_0}$, associated to the convex functional ${\cal J}_{\mathcal I_0}^p$.  As we will see, this implies that  the nonconvex functional ${\cal J}_r^p$ behaves \emph{locally} like a convex functional in neighborhoods of fixed points $u = \mathbb{H}(u)$, including the global minimizers of ${\cal J}_r^p$.

\subsection{On the nonexpansiveness and convergence for $T$ injective}
Given that $\mathbb{H} ({ u}^n) = \mathbb{U}_{\mathcal I_0} ({ u}^n)$ after a finite number of iterations, we can use well-known tools from convex analysis to prove that the sequence $({ u}^n)_{n \in \mathbb N}$ converges.
If the operator $T^*T: \ell_2(\mathcal I) \rightarrow  \ell_2(\mathcal K) $ is invertible, or, equivalently, if the operator $T$ maps onto its range and has a trivial null space -- as, for example, does the discrete pseudoinverse $D_{h}^{\dagger}$ in the 1D Mumford-Shah approximation -- then the mapping $\mathbb{U}_{\mathcal I_0}$ has the nice property of being a contraction mapping, so that a direct application of the Banach fixed point theorem ensures exponential convergence of the sequence $( u^n)_{n \in \mathbb N}$ after fixation of the index sets. \begin{theorem} \label{contract}
Suppose $T: \ell_2(\mathcal I) \rightarrow \ell_2(\mathcal K)$ maps onto $\ell_2(\mathcal K)$ and has a trivial null space.  Let $\delta > 0$ be a lower bound on the spectrum of $T^* T$.
 Then the sequence 
\begin{equation}
{ u^{n+1}} = \mathbb{H} ({ u^n}), 
\end{equation}
as defined in $\eqref{u}$, is guaranteed to converge in norm.  In particular, after a finite number of iterations $N \in \mathbb N$, this mapping takes the form
\begin{equation}
{ u}^{N+m} =  \mathbb{U}_{\mathcal I_0}^m ({{ u}}^N), \quad m \in \mathbb N \setminus \{0\},
\end{equation} 
and the sequence $( u^n)_{n \in \mathbb N}$ converges to the unique fixed point $\bar{u}$ of the map $\mathbb{U}_{\mathcal I_0}$.  
Moreover, after fixation of of the index set $\mathcal I_0$, the rate of convergence becomes exponential:
\begin{equation}
\| { u}^{N+m} - { \bar{u}} \|_{\ell_2(\mathcal I)} =\| \mathbb{U}_{\mathcal I_0}^m ({ u}^{N}) - \mathbb{U}_{\mathcal I_0}^m ({ \bar{u}}) \|_{\ell_2(\mathcal I)}  \leq (1-\delta)^m\| {{ u}}^N - { \bar{u}} \|_{\ell_2(\mathcal I)}, \quad m \in \mathbb N \setminus \{0\}.
\end{equation} 
\end{theorem}
The proof of Theorem  \ref{contract} is deferred to the Appendix.  

\subsection{Convergence for general operators $T$}
Unfortunately, if $T^*T$ is not invertible (that is, if $\delta = 0$ belongs to its nonnegative spectrum), then the map $\mathbb{U}_{\mathcal I_0}$ is not necessarily a contraction, and we can no longer apply the Banach fixed point theorem to prove convergence of the sequence $( u^n)_{n \in \mathbb N}$.  
However, as long as $\|T\| < 1$, we observe by following the proof of Theorem \eqref{contract} that $\mathbb{U}_{\mathcal I_0}$ is still \emph{non-expansive}, meaning that for all $v, v' \in \ell_2(\mathcal I)$, $\| \mathbb{U}_{\mathcal I_0} (v) - \mathbb{U}_{\mathcal I_0} (v') \|_{\ell_2(\mathcal I)} \leq \| v - v' \|_{\ell_2(\mathcal I)}$.  The following Opial's theorem $\cite{op67}$, here reported adjusted to our notations and context, gives sufficient conditions under which non-expansive maps admit convergent successive iterations:
\begin{theorem}[Opial's Theorem]
Let the mapping $\mathbb{A}$ from $\ell_2(\mathcal I)$ to $\ell_2(\mathcal I)$ satisfy the following conditions:
\begin{enumerate}
\item $\mathbb{A}$ is asymptotically regular: for all $v \in \ell_2(\mathcal I)$, $\| \mathbb{A}^{n+1} (v) - \mathbb{A}^n (v) \|_{\ell_2(\mathcal I)} \rightarrow  0$ for $n \to \infty$;
\item $\mathbb{A}$ is non-expansive: for all $v, v' \in \ell_2(\mathcal I)$, $\| \mathbb{A} (v) - \mathbb{A} (v') \|_{\ell_2(\mathcal I)} \leq \| v - v' \|_{\ell_2(\mathcal I)}$;
\item the set $\operatorname{Fix}(\mathbb A)$ of the fixed points of $\mathbb{A}$ in $\ell_2(\mathcal I)$ is not empty.
\end{enumerate}
Then, for all $v \in \ell_2(\mathcal I)$, the sequence $( \mathbb{A}^n (v))_{n \in \mathbb{N}}$ converges weakly to a fixed point in $\operatorname{Fix}(\mathbb A)$. 
\label{opial}
\end{theorem}

\noindent In fact, we already know that $\mathbb{U}_{\mathcal I_0}$ is asymptotically regular, in addition to being nonexpansive -  this follows by application of Lemma $\ref{l1}$ and Lemma $\ref{l2}$ to the functional ${\cal J}_{\mathcal I_0}^p$.  Thus, in order to apply Opial's theorem, it remains only to show that $\mathbb{U}_{\mathcal I_0}$ has a fixed point; that is, that there exists a point $\bar{u} \in \ell_2(\mathcal I)$ for which $${ \bar{u}} = \mathbb{U}_{\mathcal I_0}({ \bar{u}}). $$
In more detail, we must prove the existence of a vector ${ \bar{u}} \in \ell_2(\mathcal I)$ satisfying
\begin{eqnarray}
\bar{u}_i &=& \left\{
\begin{array}{ll} F_p^{-1}([(I - T^*T)\bar u + T^*g]_i), & \textrm{if } i \in \mathcal I_0 \\  
\big( (I - T^*T)\bar{u} + T^*g \big)_i, & \textrm{if } i \in \mathcal I_1
\end{array} \right.
\label{seperate2}
\end{eqnarray}
The following lemma gives a simple yet useful characterization of points satisfying the fixed point relation $\eqref{seperate2}$:
\begin{lemma}
Suppose $p > 1$.  A vector ${ \bar{u}} \in \ell_2(\mathcal I)$ satisfies the fixed point relation ${ \bar{u}} = \mathbb{U}_{\mathcal I_0}({ \bar{u}})$ if and only if
\begin{equation}
\label{wha}
 \big[ T^*(g - T\bar{u}) \big]_i  = \left\{
\begin{array}{ll}
0, &  i \in \mathcal I_1 \\  
F_p(\bar{u}_i) - \bar{u}_i, & i \in \mathcal I_0,
\end{array}
\right.
\end{equation}
Alternatively, if $p = 1$ and $r \geq 1/4$,  ${ \bar{u}} = \mathbb{U}_{\mathcal I_0}({ \bar{u}})$ is satisfied if and only if
\begin{equation}
   \left\{
\begin{array}{ll}\label{wah}
\big[ T^*(g - T\bar{u}) \big]_i \in [-1/2, 1/2] , & i \in \mathcal I^a_0, \\  
\big[ T^*(g - T\bar{u}) \big]_i = 1/2\sgn{\bar{u}_i}, & i \in \mathcal I^b_0, \\
\big[ T^*(g - T\bar{u}) \big]_i =0, &  i \in \mathcal I_1,
\end{array}
\right.
\end{equation}
where in $\eqref{wah}$, the index set $\mathcal I_0$ is split into
\begin{itemize}
\item $\mathcal I^a_0 =  \left \{ i \in \mathcal I_0: |\bar{u}_i| \leq 1/2 \right \}$, and
\item $\mathcal I^b_0 =  \left \{ i \in \mathcal I_0: 1/2 < |\bar{u}_i| \leq r + 1/4 \right \}$.
\end{itemize}
\label{l4}
\end{lemma}
Again, recall the notation $ F_p(t) = t + \frac{p}{2}\sgn{t}{|t|}^{p-1}$, and observe that the fixed point relation \eqref{wha} has a very simple expression when $p = 2$.   The proof of Lemma \ref{l4} is given in the Appendix.
\\
\\
The fixed point characterization of Lemma $\ref{l4}$ will be crucial in the following theorem that ensures the existence of a fixed point  $\bar{{ u}} = \mathbb{U}_{\mathcal I_0} (\bar{{ u}})$.  We remind the reader that until now, all of the results of Section $1.3$ remain valid in the infinite-dimensional setting $| {\cal I} | = \infty$.  From this point on, however, certain results will only hold in finite dimensions; for clarity, we will account each such situation explicitly. 

\begin{prop}
In finite dimensions $| {\cal I} | < \infty$, then there exist (global) minimizers of the convex functional,
\begin{equation}
 {\cal{J}}_{\mathcal I_0}^p({ u}) =  \| T{ u} - { g} \|_{\ell_2(\mathcal K)}^2 +  \sum_{i \in \mathcal I_0} |u_i|^p,
 \label{lp} 
\end{equation}
for all $p \geq 1$, and any minimizer $\bar{{ u}}$ of ${\cal{J}}_{\mathcal I_0}^p$ satisfies the fixed point relation $\bar{{ u}} = \mathbb{U}_{\mathcal I_0}(\bar{{ u}})$.  Restricted to the range $1 \leq p \leq 2$, the statement is true also in the limit $| {\cal I} | = \infty$.  
\label{l5}
\end{prop}
\begin{proof}
In the finite-dimensional setting, minimizers necessarily exist for all $p \geq 1$ according to Proposition \ref{firstlemma}.   We now consider the general case. Consider the unique decomposition $u = u_0 + u_1$ into a vector $u_0$ supported on $\mathcal I_0$ and another $u_1$ supported on $\mathcal I_1$, i.e., the vectors $u_0 \in \ell_2^{\mathcal I_0}(\mathcal I) := \{u \in \ell_2(\mathcal I): u_i=0, \quad i \in \mathcal I_1\}$ and $u_1 \in \ell_2^{\mathcal I_1}(\mathcal I) := \{u \in \ell_2(\mathcal I): u_i=0, \quad i \in \mathcal I_0\}$. 
Let $\mathcal P: u \rightarrow u_1$ and $\mathcal P^{\perp} = \mathcal I - \mathcal P: u \rightarrow u_0$ denote the orthogonal projections onto the subspaces  $\ell_2^{\mathcal I_1}(\mathcal I)$ and $\ell_2^{\mathcal I_0}(\mathcal I)$, respectively.   Consider the operators $T_0 = T \mathcal P^{\perp}$ and $T_1 = T \mathcal P$; note that clearly $T = T_0 + T_1$ is satisfied.
%Write $T_0: {\cal H}_0 \rightarrow {\cal{\tilde{H}}}_0$ and $T_1: {\cal H}_1 \rightarrow {\cal{\tilde{H}}}_1$.  
\mnew{ The functional \eqref{lp} can be re-written with this decomposition according to
\begin{equation}
 {\cal{J}}_{\mathcal I_0}^p({ u_0} + { u_1}) =  \| T_0 u_0 + T_1 u_1 - g \|_{\ell_2(\mathcal K)}^2 +  \|u_0\|_{\ell_p^{\mathcal I_0}(\mathcal I)}^p
\end{equation}  
where $\|z \|_{\ell_p^{\mathcal I_0}(\mathcal I)}:= \left (\sum_{i \in \mathcal I_0} |z_i|^p\right)^{1/p} $ is the $\ell_p$-norm on vectors supported on $\mathcal I_0$.
\\
\\
\noindent Let $\mathcal P_1$ be the orthogonal projection onto the range of $T_1$ in $\ell_2(\mathcal K)$ (not to be confused with $\mathcal P$, which operates on the space $\ell_2(\mathcal I)$}) and let $\mathcal P_1^{\perp} = \mathcal I  - \mathcal P_1$ be the orthogonal projection in $\ell_2(\mathcal K)$ onto the orthogonal complement of the range of $T_1$.  Then, fixing $u_0 \in  \ell_2^{\mathcal I_1}(\mathcal I)$, the vector $\mathcal P_1(g - T_0 u_0) \in \operatorname{range}(T_1) \subset \ell_2(\mathcal K)$ is the solution to the minimization problem
\begin{equation}
\mathcal P_1(g - T_0 u_0)= \arg \min_{v \in \operatorname{range}(T_1)} \| v - (g - T_0 u_0) \|^2_{\ell_2(\mathcal K)}, 
\end{equation}
so that minimizers of the functional ${\cal{F}}: \ell_2^{\mathcal I_0}(\mathcal I) \rightarrow \mathbb{R}^+$ defined by
\begin{eqnarray}
 {\cal{F}}(v) &=&  \| T_0 v + \mathcal P_1(g - T_0 v) - g \|^2_{\ell_2(\mathcal K)} + \|v \|_{\ell_p^{\mathcal I_0}(\mathcal I)}^p \nonumber \\
 &=& \| K v - y \|^2_{\ell_2(\mathcal K)} + \|v \|_{\ell_p^{\mathcal I_0}(\mathcal I)}^p
 \label{DDD}
\end{eqnarray}
with $K := \mathcal P_1^{\perp} T_0$, and $y := \mathcal P_1^{\perp} g$, will yield minimizers of ${\cal{J}}_{\mathcal I_0}^p$.   Functionals of the form $\eqref{DDD}$ were studied in $\cite{DDD}$; there, it is shown that as long as $1 \leq p \leq 2$, ${\cal{F}}$ has minimizers, and any minimizer $\bar{v}$ can be characterized by the fixed point relation
\begin{equation}
\bar{v}_i = F_p^{-1}([(I - K^*K)\bar{v} + K^*y]_i), \quad i \in \mathcal I_0;
\label{usual}
\end{equation}
(recall that $F_p^{-1}$ is the inverse of the function $F_p(t) = t + \frac{p}{2}\sgn t |t|^{p-1}$).
\\
 In the finite-dimensional setting $| \cal I | < \infty$, the Euler-Lagrange equations corresponding to minimizers of the convex functional ${\cal{F}}$ as in \eqref{DDD} imply the same fixed point relation \eqref{usual} also, for all $p \geq 1$. 
\\
By Lemma $\ref{l4}$, the characterization \eqref{usual} is equivalent to the condition
\begin{itemize}
\item $p > 1$:
\begin{equation}
 \big[ K^*(y - K\bar{v}) \big]_i  =  \frac{p}{2}\sgn{\bar{v}_i} |\bar{v}_i|^{p-1}, 
 \label{ey}
\end{equation}
\item $p = 1$: 
\begin{equation}
 \left\{
\begin{array}{ll}
 \big[ K^*(y - K\bar{v}) \big]_i  \in [-1/2, 1/2], & \textrm{if } |\bar{v}_i| \leq 1/2,  \\  
 \big[ K^*(y - K\bar{v}) \big]_i  = 1/2\sgn{\bar{v}_j}, & \textrm{if } 1/2 < |\bar{v}_i| \leq r + 1/4.
\end{array}
\right., \quad i \in \mathcal I_0.
\label{ew}
\end{equation}
\end{itemize}
Making the identification $\bar{u}_0 = \bar{v}$ and $T_1\bar{u}_1 = \mathcal P_1(g - T_0 \bar{v})$, and rewriting $K = \mathcal P_1^{\perp} T_0$, and $y = \mathcal P_1^{\perp} g$, the relations $\eqref{ey}$ and $\eqref{ew}$ imply the full fixed point characterization in Lemma \ref{l4}. 
\end{proof}

\begin{remark}
\emph{The restriction $p \leq 2$ that is necessary for the results of this paper in the infinite dimensional setting $| \mathcal I | = \infty$ was only used in the proof of Theorem $\ref{l5}$, where it comes from $\cite{DDD}$ and is needed there to prove the existence of minimizers of functionals ${\cal{F}}$ of the form $\eqref{DDD}$.   If that proof can be extended to functionals of the form $\eqref{DDD}$ for general $p \geq 1$, then the restriction $p \leq 2$ can be dropped in the current paper.  
For instance, if we additionally require that $T$ is a bounded operator from $\ell_p(\mathcal I)$ to $\ell_2(\mathcal I)$ for $1 \leq p < \infty$ then the existence of minimizers would be guaranteed also for $1 \leq p < \infty$ and $| \mathcal I | = \infty$. In this case we could consider a minimizing sequence $(v^k)$ of $\mathcal F$, which is necessarily bounded in $\ell_p$. Therefore, there exists a subsequence $(v^{k_h})$ which weakly converges in $\ell_p$ to a point $v^*$. This also implies the weak convergence of the sequence $K v^{k_h}$ in $\ell_2$; note that $\langle K v^{k_h}, w \rangle_{\ell_2 \times \ell_2} = \langle v^{k_h}, K^* w \rangle_{\ell_p \times \ell_{p'}}$, for $1/p+1/p' =1$. By Fatou's lemma we obtain $\mathcal F (v^*) \leq \lim \inf_h \mathcal F (v^{k_h})$ and $v^*$ is a minimizer of $\mathcal F$.
However, we still require that $p \geq 1$ for the proof of Proposition $\ref{thmthrs}$ and for the results of the next section to hold.  }
\end{remark}

\noindent Combining the results from this section, we obtain:

\begin{theorem}\label{mainth}
Suppose $1 \leq p \leq 2$.  Starting from any ${ u}^0$ satisfying ${\cal{J}}^p_r({ u}^0) < \infty$, the sequence $({ u}^n)_{n \in \mathbb N}$ defined by ${ u}^{n+1} = \mathbb{H}^n ({ u}^0)$ as in \eqref{u} will converge weakly to a vector $\bar{{ u}} \in \ell_2(\mathcal I)$ that satisfies the fixed point condition,
\begin{enumerate}
\item $|\bar{u}_i| \geq \lambda'(r,p)$, if $i \in \mathcal I_1 = \{j \in \mathcal I: |\bar{u}_j| > r\}$
\item $|\bar{u}_i| \leq F_p^{-1}(\lambda'(r,p))$, { for $p>1$,} if $i \in \mathcal I_0 = \{j \in \mathcal I: |\bar{u}_j| \leq r\}$, and
\item 
\begin{enumerate}
\item If $p > 1$: 
\begin{equation}
\label{wha2}
 \big[ T^*(g - T\bar{u}) \big]_i  = \left\{
\begin{array}{ll}
0, & \textrm{ if }  |\bar{u}_i| \geq  \lambda'(r,p) \\
F_p(\bar{u}_i) - \bar{u}_i, & \textrm{ if } |\bar{u}_i| \leq  \lambda'(r,p) - \delta(r,p)
\end{array}
\right.
\end{equation}
\item If $p = 1$ and $r \geq 1/4$: 

\begin{equation}
  \left\{
\begin{array}{ll}
\big[ T^*(g - T\bar{u}) \big]_i  \in [-1/2, 1/2] , & |\bar{u}_i| \leq 1/2 \\  
\big[ T^*(g - T\bar{u}) \big]_i  = 1/2\sgn{\bar{u}_i}, & 1/2 < |\bar{u}_i| \leq r - 1/4. \\
\big[ T^*(g - T\bar{u}) \big]_i  = 0, &  |\bar{u}_i| > r + 1/4.
\label{wah2}
\end{array}
\right.
\end{equation}
\end{enumerate}
\end{enumerate}
If the index set $| \cal I | < \infty$ is finite dimensional, the theorem holds for all $p \geq 1$.  
\label{bigthm2}
\end{theorem}
\begin{proof}
By Lemma $\ref{l3}$, the map $u^{n+1} = \mathbb{H}(u^n)$ becomes equivalent to a map of the form
$u^{n+1}  = \mathbb{U}_{\mathcal I_0} (u^n)$ after a finite number of iterations $N \in \mathbb N$.  By Lemma $\ref{l3}$ and Proposition $\ref{thmthrs}$, the subset $\mathcal I_0 \subset \mathcal I$ separates $\mathcal I$ in the sense that, for all $n \geq N$, 
\begin{itemize}
\item $|u_i^n| < F_p^{-1}(\lambda'(r,p))$, if $i \in \mathcal I_0$,
\item $|u_i^n| > \lambda'(r,p)$, if $i \in \mathcal I_1 = \mathcal I \setminus \mathcal I_0$.
\end{itemize}
That the sequence $(u^n)_{n \in \mathbb N}$ converges to a fixed point of the map $\mathbb{U}_{\mathcal I_0} $ follows from Opial's theorem applied to the map $\mathbb{U}_{\mathcal I_0}$:   
\begin{enumerate}
\item the asymptotic regularity of $\mathbb{U}_{\mathcal I_0}$ is a consequence of Lemmas $\ref{l1}$ and $\ref{l2}$;
\item the nonexpansiveness of $\mathbb{U}_{\mathcal I_0}$ follows from the proof of Theorem \eqref{contract}, and
\item Theorem $\ref{l5}$ guarantees that the set of fixed points of $\mathbb{U}_{\mathcal I_0}$ in $\ell_2(\mathcal{I})$ is nonempty.
\end{enumerate}
The limit $\bar{u}$ of the sequence $(u^n)$ will satisfy the fixed point conditions of Lemma \ref{l4}.  Since weak convergence implies component-wise convergence, it follows for all $i \in \mathcal I_0$ that
\begin{eqnarray}
|\bar{u}_i| &=& \lim_{n \rightarrow \infty} |u^n_i| \nonumber \\
&\leq& \lambda'(r,p) - \delta(r,p)
\end{eqnarray}
and the respective lower bound $|u_i^n| \geq \lambda'(r,p)$ holds analogously for $i \in \mathcal I_1$.
\end{proof}

\section{On minimizers of ${\cal J}^p_r$}

We are now in a position to explore the relationship between limit vectors $\bar{u}$ of the iterative thresholding algorithm \eqref{u} and minimizers of the free-discontinuity functional  ${{\cal{J}}}_r^p$ \eqref{generalform}.   As a first but important result in this direction, 

\begin{theorem}
A point $\bar{u}$ satisfying the fixed point relation of Theorem $\ref{bigthm2}$ is a local minimizer of the functional ${\cal{J}}^p_r$ defined in $\eqref{generalform}$. 
\label{mintheorem} 
\end{theorem}
The proof of Theorem \ref{mintheorem} is omitted at present but can be found in the Appendix.   This result should not be surprising, however.   Due to the separation of the entries of any fixed point $\bar u$, such that $\bar u_i < r < \bar u_j$ for $i \in \mathcal I_0$ and $j \in \mathcal I_1$,  we have also $\mathcal I_0 \equiv \{i \in \mathcal I: |u_i| \leq r\}$ and $\mathcal I_1 \equiv \{j \in \mathcal I: |u_j| > r\}$ for all $u \in B(\bar u,\varepsilon(r))$, where $B(\bar u,\varepsilon(r))$ is a ball around an equilibrium point $\bar u$ of radius $\varepsilon(r)>0$ sufficiently small.  On this neighborhood $B(\bar u,\varepsilon(r))$ of $\bar{u}$, the functional $\mathcal J_r^p$ is convex.  Since $\bar{u}$ is obtained as the limit of a sequence $(u^n)$ in $B(\bar u,\varepsilon(r))$ for which the sequence ${\cal J}^p_r(u^n)$ is nonincreasing, one would expect that $\bar{u}$ minimizes ${\cal J}^p_r(u^n)$ within this neighborhood.  
\\
\\
More surprising is that global minimizers of ${\cal J}^p_r$ are also fixed points, as shown in the following theorem.   Even though the existence of such minimizers is only guaranteed in the finite-dimensional setting (see Proposition \ref{firstlemma}), the following result is not restricted as such.

\begin{theorem}[Global minimizers of ${\cal J}_r^p$ are fixed points $\bar{u} = \mathbb{H}(\bar{u})$]
{ Any global minimizer $u^*$ of ${\cal{J}}^p_r$ satisfies the fixed point condition of the map $\mathbb{H}$ that is given in Theorem \ref{bigthm2}. } 
\label{globalmin}
\end{theorem}

The proof of Theorem \ref{globalmin} is rather long and we defer it to the Appendix.    We reiterate once more that on a ball $B(\bar u,\varepsilon(r))$ around an equilibrium point $\bar u$ of radius $\varepsilon(r)>0$ sufficiently small, the functional $\mathcal J_r^p$ is convex; following the proof of Theorem \ref{globalmin}, we see that $\mathcal J_r^p$ is in fact \emph{strictly} convex whenever $\bar u=u^*$ is a global minimizer, since the restriction of $T$ to the subspace $\ell_2^{\mathcal I_1}(\mathcal I) \subset \ell_2(\mathcal I)$ of vectors with support in $\mathcal I_1$ must be an injective operator in this case. Hence a global minimizer is necessarily an isolated minimizer, whereas we cannot ensure the same property for local minimizers if $T$ has a nontrivial null-space; in this case, local minimizers may form continuous sets, as it is shown in the bottom-right box of Figure \ref{patterns}. We conclude the following remark.
\begin{corollary}\label{isolated}
Minimizers of $\mathcal J_r^p$ are isolated.  
\end{corollary}  
%
%Corollary \eqref{isolated} marks an important difference in nature between energy functionals $\mathcal J =  \| T u - g \|_{\ell_2(\mathcal K)} + R(u)$ having free-discontinuity type regularization $R(u) = \sum_{j \in \mathcal I} \min\{ | u_j|^p, r^p \}$, and $\ell_1$-regularized functionals (that is, where $R(u) = \| u \|_{\ell_1(\mathcal I)}$).  In the latter situation, minimizers may form continuous sets \cite{DDD}.  

\mnew{
\section{2-D free-discontinuity inverse problems and a projected gradient method}
As presented in Subsection \ref{2dcase}, the minimization of the discrete functionals for 2-D free-discontinuity inverse problems has the general form
\begin{eqnarray}
&&
 \left\{ \begin{array}{llll}
\textrm{Minimize} &  {\cal{J}}_{r}^p(u) :=  \| T u - g\|_{\ell_2(\mathcal K)}^2 + \sum_{i \in \mathcal I} \min \left \{|u_{i}|^p, r^p  \right \} 
\\
 \textrm{subject to} & {\cal {Q}} u = 0,\\
\end{array}
\right.
\label{MS2dgeneral}
\end{eqnarray}
where ${\cal Q}:\ell_2(\mathcal I) \to \ell_2(\mathcal K')$ is a suitable bounded linear operator. 
\\
\\
\noindent We can not directly generalize the analysis of the previous sections to $\eqref{MS2dgeneral}$, as the introduction of surrogate functionals does not decouple the constraint ${\cal Q}u = 0$.  
However, when the index set $\mathcal I$ is finite dimensional, we can still say something.  For ease of presentation, we will assume $p = 2$ throughout this section.  
\\
\\
\noindent First, recall that the partition argument of Theorem $\ref{existmin}$ guarantees that the constrained minimization problem \eqref{MS2dgeneral} has a minimizer.  Again, one could in theory obtain such a minimizer by computing a minimizer $u({\cal I}_0)$ for each subset ${\cal I}_0 \subset {\cal I} = \{1,2,..., N\}$.  Of course, such an algorithm is computationally infeasible as the number of subsets of the index set $\{1,2, ..., N \}$ grows exponentially with the dimension $N$ of the underlying space. 
\\
\\
\noindent We propose instead the following more practical projected gradient algorithm: for any initial $u^0$, iterate
\begin{equation}
\label{2DMSalg}
u^{n+1} = \mathcal P_{\ker(\mathcal Q)} \left [ H_{(2,r)}(u^n + T^*(g - T u^n)) \right ],
\end{equation}
where $\mathcal P_{\ker(\mathcal Q)}$ is the orthogonal projection onto the null-space of $\mathcal Q$. This projection can be easily computed explicitly by
\begin{eqnarray*}
\mathcal P_{\ker \mathcal Q} &=& I - \mathcal Q^\dagger \mathcal Q\\
&=& I - \mathcal Q^* (\mathcal Q \mathcal Q^*)^{-1} \mathcal Q,
\end{eqnarray*}
where the latter equality holds whenever $\mathcal Q$ is a full-rank matrix, as the one associated to the Schwartz conditions \eqref{comp}.  The analysis of the algorithm \eqref{2DMSalg} is beyond the scope of this paper; nevertheless, note that locally around any minimizer, the functional $\mathcal J_r^2$ is convex, and that projected gradient iterations are well-known methods for constrained minimization of (non-smooth) convex functionals, see for instance \cite{aliuso98}.   

\section{Numerical Experiments}

\subsection{Dynamical systems, stability, and equilibria}
Iterative thresholding algorithms have a natural interpretation as discrete-time dynamical systems with nonsmooth right-hand-side, and can be associated to continuous dynamical systems of the type:
\begin{eqnarray*}
\dot u(t) &=& F(u(t),t)\\&=&\tau \left ( H_{(p,r)}( u(t) +  T^*(g - T u(t))) - u(t) \right), \quad t \geq t_0, \quad \tau >0.
\end{eqnarray*}
The study of the existence, uniqueness, stability, and long-time behavior  of these ODE's is of fundamental interest in order to clarify also the stability properties of iterative thresholding algorithms. Indeed, other than soft-thresholding iterations \cite{DDD}, the corresponding right-hand-side is not Lipschitz continuous and can even be  discontinuous, as is the case for free-discontinuity problems. In  \cite{br88,fi88} conditions are established for the existence, uniqueness, and continuous dependence on the initial data (at finite time) { of solutions of dynamical systems with discontinuous right-hand-side}. However, very little is known about long-time properties of such dynamical systems and about the nature of their equilibrium points. \\
\begin{figure}
\[\begin{array}{cc}
  \centering
  \includegraphics[width=7cm]{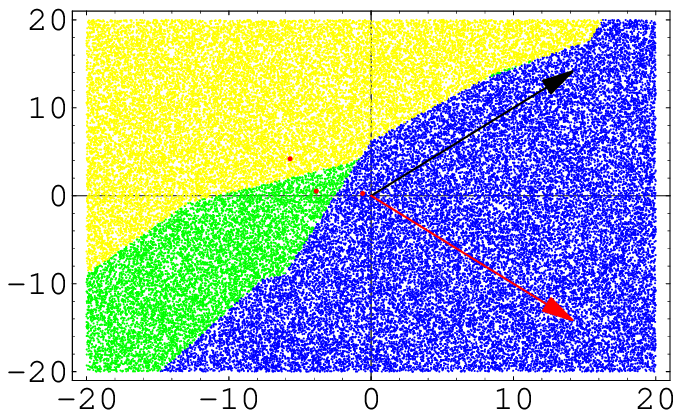} &  
  \includegraphics[width=7cm]{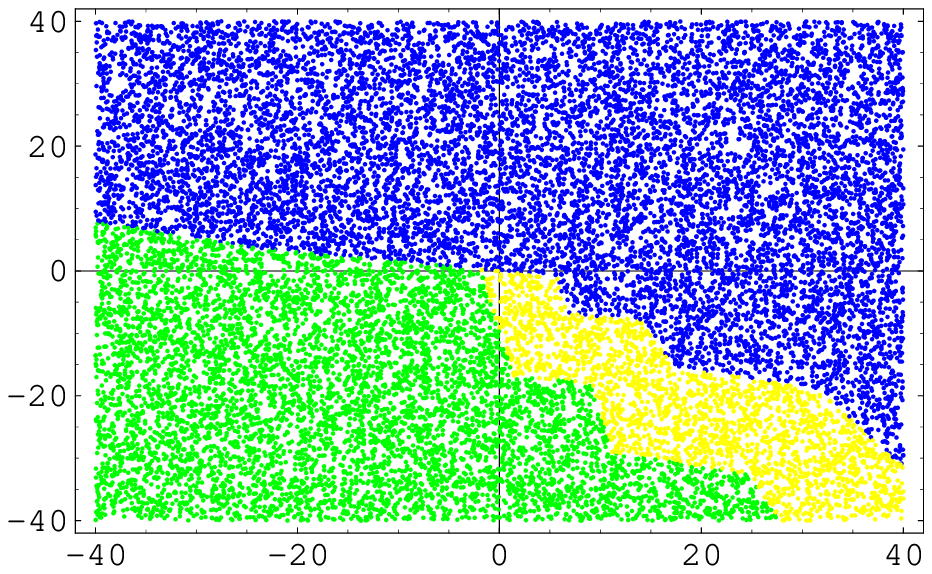} \\
  \includegraphics[width=7cm]{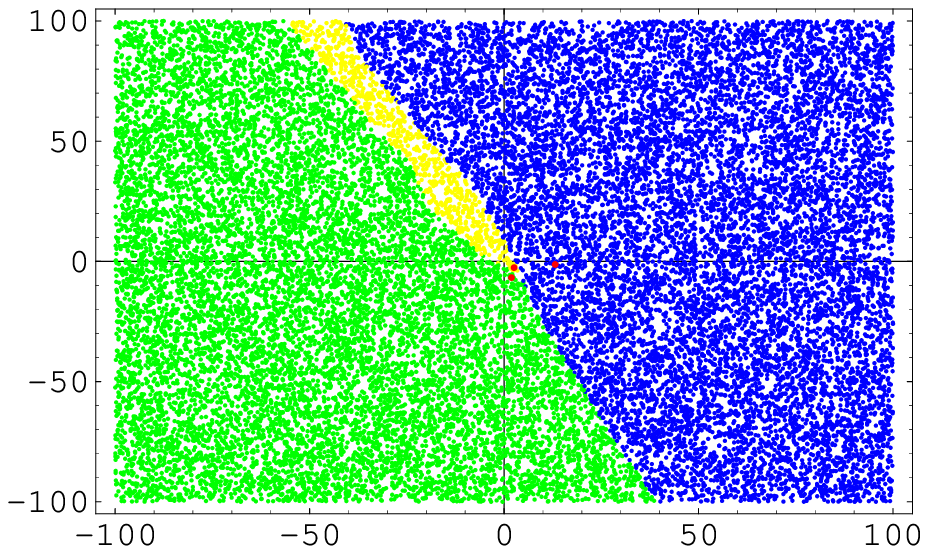} &
  \includegraphics[width=7cm]{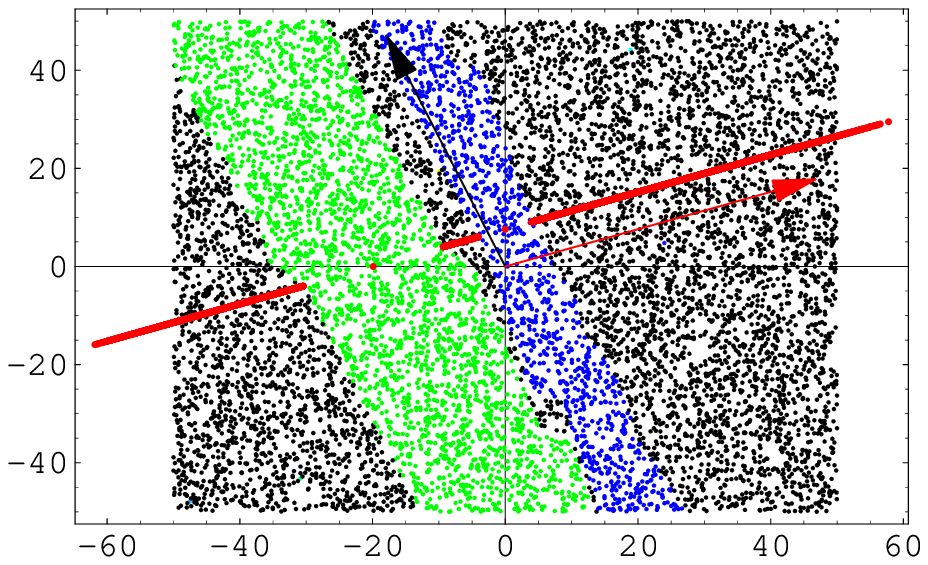} 
\end{array}\]
\caption{We show patterns in $\mathbb R^2$ formed by initial points $u^0$ colored according to the corresponding equilibria computed as limits of the iterative thresholding algorithm \eqref{map3}. For invertible $2\times 2$ squared matrices $T$, the equilibria are isolated and the region of initial points for which \eqref{map3} converges to a given equilibrium point do partition the space into sets which might be disconnected. Structures of the partition generated by different matrices $T$ are exemplified in the top boxes and in the bottom-left one. In the bottom-right box we show the pattern related to iterations where the  $2\times 2$ squared matrix $T$ has nontrivial null-space. We can see again that global minimizer are isolated and correspond to the points on the axes, whereas local minimizers are continuously distributed along an affine space generated by the kernel of $T$. It is not difficult to show that this structure always occurs for such matrices.} 
\label{patterns}
\end{figure}

\noindent For several continuous thresholding functions, such as the ones introduced in \cite{DDD,fora06,fora07}, one can easily show, for instance by means of $\Gamma$-convergence arguments, that equilibrium points depend continuously on the parameters of the thresholding, see, e.g., \cite[Theorem 5.1]{fora07}. Nevertheless, for discontinuous thresholding functions $H_{(p,r)}$ such as those studied in this paper, sudden bifurcation phenomena and instabilities do appear in general. Figure \ref{patterns} shows that multiple equilibrium points can exist for these thresholding operators and their number may depend discontinuously on the thresholding shape parameters. Moreover, as established in Theorem \ref{globalmin}, global minimizers of $\mathcal J_r^p$ are always stable equilibria and isolated points, while local minimizers can be unstable equilibria and form a continuous set, as shown  in the bottom-right box of Figure \ref{patterns}.

\subsection{Denoising and segmentation of 1-D signals and digital images}
\label{denoising}

In this subsection, we are concerned with numerical experiments in the use of an iterative thresholding algorithm for the minimization of
\begin{equation}
{\cal{J}}_{r,\gamma}^2 ({u}) :=  \| D_h^\dagger u  - g \|^2_{\ell_2} +  \gamma \sum_{i=1}^N \min \{ u_i^2, r^2 \},
\end{equation}
modelling problems of denoising and segmentation.
\begin{figure}[htp]
\begin{center}
\includegraphics[width=6in]{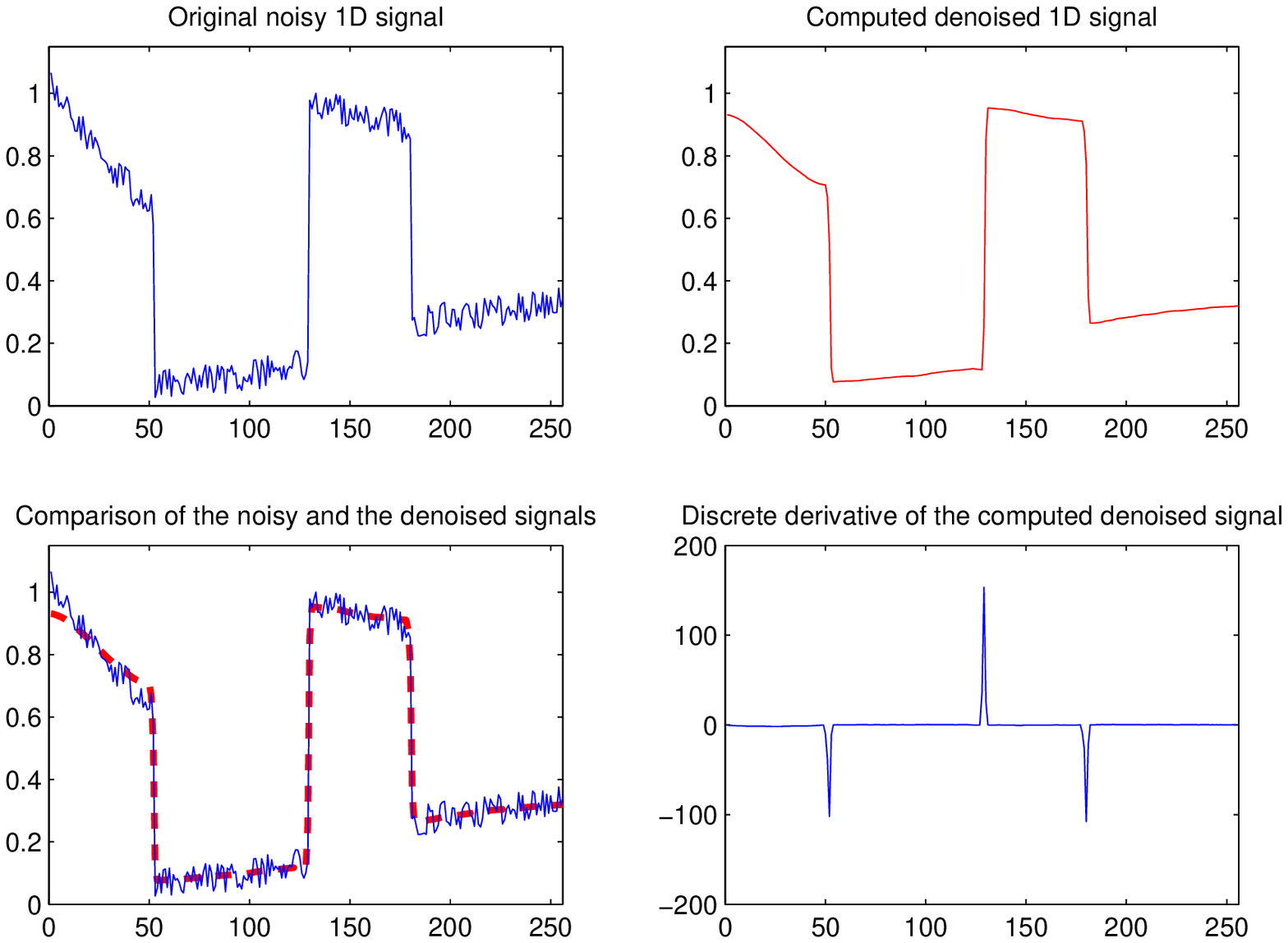}
\caption{We show the application of the  iterative thresholding algorithm \eqref{map3} for the classical denoising problem of 1-D signals where $K=I$ in \eqref{regMS}, and hence $T = D_h^\dagger$. The thresholding parameters used for the numerics are $r=2.2$ and $\gamma=0.002$.}\label{1Dsignal}
\end{center}
\end{figure}

\begin{figure}[htp]
\begin{center}
\includegraphics[width=6in]{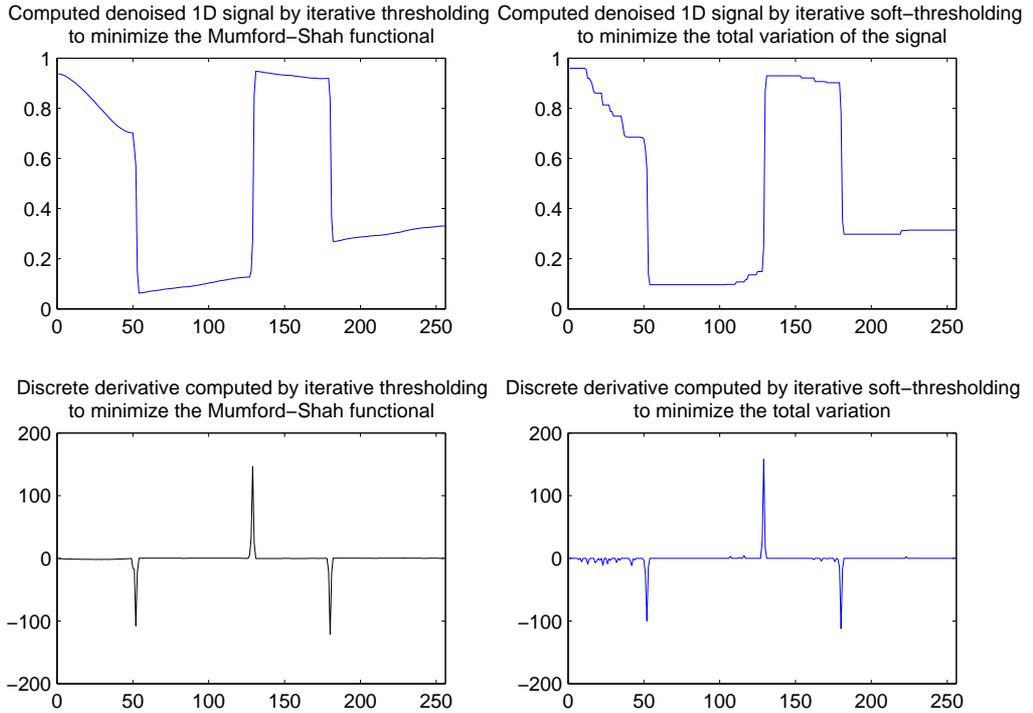}
\caption{A comparison of the denoising of the signal in Figure \ref{1Dsignal} by means of the algorithm \eqref{map3} and by iterative soft-thresholding \cite{DDD} applied to discrete derivatives. We can appreciate how the algorithm \eqref{map3} promotes piecewise smooth solutions, whereas the iterative soft-thresholding promotes the total variation minimization with the introduction of a `staircase effect'.  The thresholding parameters used for the numerics are $r=2.2$ and $\gamma=0.002$ for \eqref{modH}, and $\gamma=0.002$ for the soft-thresholding  \eqref{softthrs}.}\label{MS-TV}
\end{center}
\end{figure}

Note that we introduced an additional regularization parameter $\gamma>0$ which has the sole effect of modifying the thresholding function $H_{(2,r,\gamma)}$ as follows
\begin{equation}
\label{modH}
H_{(2,r,\gamma)}(z) = \left \{ 
\begin{array}{ll}
\frac{1}{1+\gamma} z, & |z| \leq \frac{r}{\sqrt{\frac{2}{(1+\gamma)^2} + 1 - \frac{2}{1+\gamma}}}\\
z, & \mbox{ otherwise.}
\end{array}
\right.
\end{equation}
This thresholding function can be again easily computed by means of an argument similar to the proof of Proposition \ref{thmthrs}. In Figure \ref{1Dsignal} and Figure \ref{2Dsignal} we show the results of applications of the iterative thresholding algorithm \eqref{map3} and the  projected gradient algorithm \eqref{2DMSalg} respectively. { In Figure \ref{MS-TV} we show a comparison of the use of the thresholding $H_{(2,r,\gamma)}$ and the soft-thresholding $S_\gamma$ (see its definition in \eqref{softthrs}); the former promotes the minimization of the Mumford-Shah constraint $MS$ and piecewise smooth solutions, whereas the latter promotes the minimization of a total variation constraint \cite{ROF}, which is also well-known to produce (almost) piecewise constant solutions with a perhaps unwanted `staircase effect'; see also \cite[Section 4]{chli97} for details.} 
\begin{figure}[htp]
\begin{center}
\includegraphics[width=6.5in]{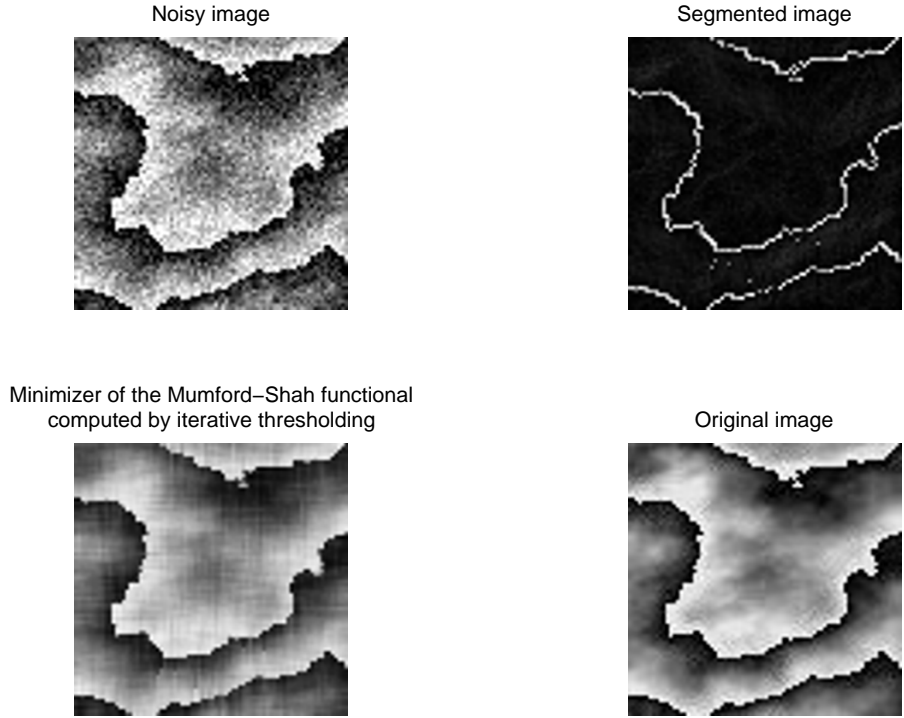}
\caption{We show the application of the  projected gradient algorithm \eqref{2DMSalg} for the classical denoising problem of digital images where $K=I$ in \eqref{regMS}, and hence $T = D_h^\dagger$.  The thresholding parameters used for the numerics are $r=5$ and $\gamma=0.005$, and the image size is $80 \times 80$. The anisotropic effects of \eqref{MS2d} are clearly visible, suggesting that for more effective image denoising, iterative thresholding on an isotropic (or direction-independent) variant of the 2D Mumford-Shah functional should be studied; see \cite{chdm99,boch00}}\label{2Dsignal}
\end{center}
\end{figure}
\subsection{Inverse problems}

As already mentioned in Subsection \ref{MS4invprob} the Mumford-Shah term  $MS ( u ) = \int_{\Omega \setminus S_u} | \nabla u |^2  +\beta \mathcal{H}^{d-1}(S_u)$ is also used for regularizing inverse problems involving operators $T$ which are not boundedly invertible. In this section we present two numerical experiments on the use of algorithms \eqref{map3} and \eqref{2DMSalg} for 1D interpolation (Figure \ref{1Dinterpolation}) and for 2D inpainting (Figure \ref{2Dinpainting})  respectively. In this case the operator $T$ is a multiplier by a characteristic function of a subdomain, i.e., $T u := \chi_D \cdot u$, for $D \subset \Omega$; see \cite{essh02} for other numerical examples previously obtained with the Mumford-Shah regularization.

\begin{figure}[htp]
\begin{center}
\includegraphics[width=6in]{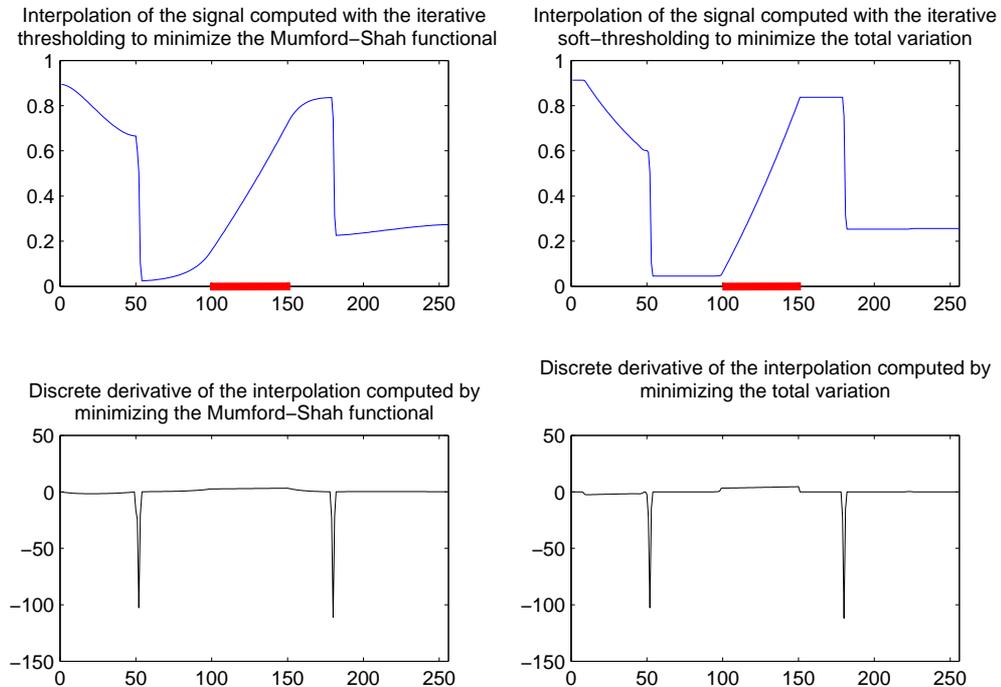}
\caption{Interpolation of an incomplete signal by means of the Mumford-Shah regularization and the total variation minimization provided by respective iterative thresholding algorithms. The red interval is the region where no information on the original signal is provided. The thresholding parameters used for the numerics are $r=2.2$ and $\gamma=0.002$  for \eqref{modH}, and $\gamma
=0.002$ for the soft-thresholding \eqref{softthrs}.}\label{1Dinterpolation}
\end{center}
\end{figure}

In Figure \ref{1Dinterpolation}  we show the reconstruction of the noiseless signal of Figure \ref{1Dsignal} provided information only out of the interval $[100,150]$ which has to be restored. On the left boxes we show the results due to algorithm \eqref{map3} and on the left ones the solution computed by iterative soft-thresholding. In the former the solution is again piecewise smooth and in the latter a (almost) piecewise constant solution is instead produced.

\begin{figure}[htp]
\begin{center}
\includegraphics[width=6in]{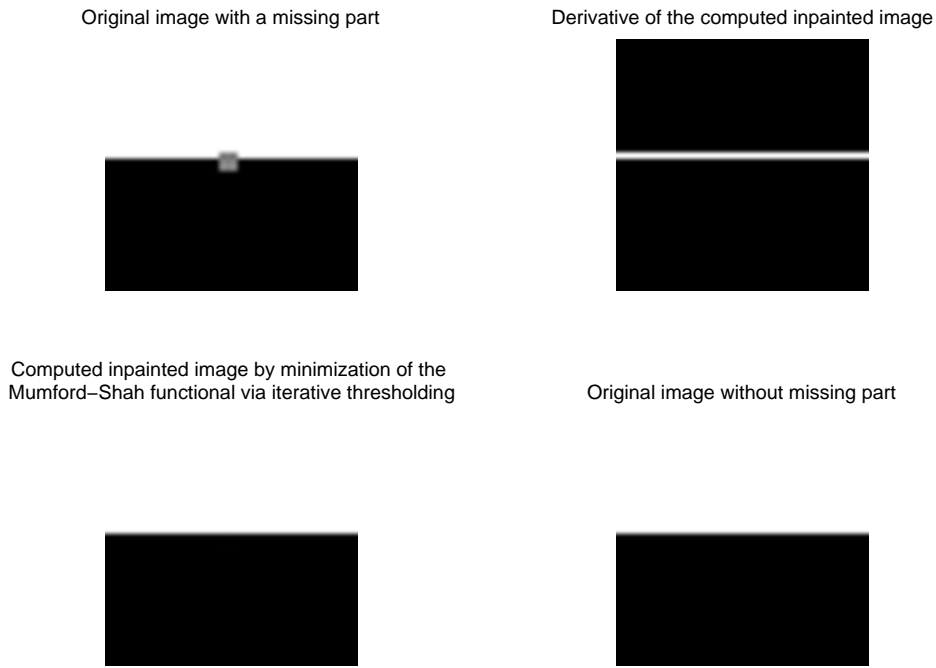}
\caption{Inpainting of a binary image by means of algorithm \eqref{2DMSalg}. The occluded discontinuity is correctly recovered as already observed in \cite{essh02}.  The thresholding parameters used for the numerics are $r=8$ and $\gamma=0.0001$, and  the image size is $40 \times 40$.}\label{2Dinpainting}
\end{center}
\end{figure}

In Figure  \ref{2Dinpainting} we show the {\it inpainting} of a binary image with a missing information right at its center which is occluding precisely a discontinuity. As already shown in \cite{essh02} the inpainting process produces minimal length connections of the discontinuity set as long as the inpainting region, i.e., the missing part, is not too large.

 \section{Appendix}

\subsection{Proof of Proposition \ref{firstlemma}}
 First, we recall Weierstrass' Theorem, which is used in the proof of Proposition $\ref{firstlemma}$ below.
\begin{theorem}[Weierstrass' Theorem]
The set of minima of a convex function $f$ over a subset $X \subset \mathbb{R}^N$ is nonempty and compact if $X$ is closed, $f$ is lower semicontinuous over $X$, and the function $\tilde{f}$, given by
\begin{equation}
\tilde{f} = \left\{ \begin{array}{ll}
f(x) & \textrm{, if } x \in X, \\
\infty & otherwise,
\end{array}
\right.
\end{equation}
is \emph{coercive}, i.e., for every sequence $( x_k) \subset X$ s.t. $\| x_k \| \rightarrow \infty$, we have $\lim_{k \rightarrow \infty} f(x_k) = \infty$.
\end{theorem}
 
% For sake of ease, we repeat now the statement of Theorem \eqref{firslemma}.
% \begin{theorem}
% Suppose $A$ is an $N \times N$ positive semidefinite matrix, and suppose $b$ and $c$ are $N \times 1$ vectors.   Suppose also that $X$ is a nonempty convex polyhedral subset of $\mathbb{R}^N$.  The convex optimization problem
% \begin{eqnarray}
%  \left\{ \begin{array}{llll}
% \textrm{minimize} & F(u) = \Big[ u^t A u + b^t u + \sum_{1 \leq j \leq N} c_j |u_j|^p \Big]  \\
%  \textrm{subject to} 
% & u \in X  \\\end{array}
% \right.  
% \label{FWgen}
% \end{eqnarray}
% admits minimizers for $p \geq 1$, as long as $F$ is bounded from below.  
% \label{firstthm}
% \end{theorem}
The following two lemmas will be helpful in the proof of Proposition \ref{firstlemma}.  

\begin{lemma}
Let $F(u)$ be a convex function defined on $\mathbb{R}^N$ having the general form $F(u) = \Big[ u^t A u + b^t u + \sum_{1 \leq j \leq N} |u_j|^p \Big]$, for some $p \geq 1$.  Fix $x$ and $d$ in  
$\mathbb{R}^N$.  If $F$ is bounded above and below on the ray $\{x + td, t \geq 0\}$, then $F$ is constant on the line $x + td$.  
\label{appendlemma1}
\end{lemma}
\begin{proof}
Let $\mu(t) = F(x + td)$, and note that $\mu$ is convex because $F$ is convex. Moreover, $\mu$ has the general form $\mu(t) = P(t) + \sum_{1 \leq j \leq N} c_j \|x_j + td_j \|^p$ where $P(t)$ is a polynomial in $t$ of order at most $2$. Without loss of generality, suppose $0 \leq \mu(t) \leq 1$ for all values of $t \in \mathbb{R}^+$.  Then there exists a sequence of points $(t_n)_{n \in \mathbb N}$, $t_n \to \infty$ for $n \to \infty$, for which $\mu(t_n)$ is a convergent sequence; let us denote the limit of this sequence by $\gamma$.

\begin{enumerate} 
\item {\bf Case 1: $1 \leq p \leq 2$}.  To repeat, 
\begin{equation}
\lim_{n \rightarrow \infty} \mu(t_n) = \lim_{n \rightarrow \infty} P(t_n) + \sum_{1 \leq j \leq N} c_j \|x_j + t_n d_j \|^p = \gamma.
\end{equation}  
Since $0 = \lim_{n \rightarrow \infty} \mu(t_n) / t_n^2$, it follows that all coefficients in $\mu(t)$ of degree 2 must vanish.  In turn, then, $0 = \lim_{n \rightarrow \infty} \mu(t_n) / t_n^p $, has the implication that for each $j$, one of the coefficients $c_j$ or $d_j$ must vanish as well.   Following in the same manner, we conclude that all linear coefficients in $\mu(t)$ also vanish, leaving only the possibility that $\mu(t) \equiv \gamma$ is a constant function.  
\item {\bf Case 2: $p > 2$}: The proof in this case is identical to that of the previous case, and as such we leave the details to the reader.  
\end{enumerate}
\end{proof}

\begin{lemma}
Suppose $F$ is a convex function defined on $\mathbb{R}^N$ that is bounded from below, and has the property that if $F$ is bounded above on a ray $\{x + td, t \in \mathbb{R}^{+} \}$, then $F$ is constant on the line $x + td$.  Then if $F$ is constant on the line $x + td$, $F$ is also constant on any parallel line $y + td$.  
\label{appendlemma2}
\end{lemma}
\begin{proof}
Let $\mu(t) = F(x + td)$ which by assumption is a constant function $\mu(t) = \gamma$, and let $v(t) = F(y + td)$.  Fix $t \in \mathbb{R}^+$, and let $z$ be the point $z = x + 2(y-x)$, i.e. $y = \frac{1}{2}x + \frac{1}{2}z$.   By convexity of $F$, we have that 
\begin{equation}
F(y + td) = F \Big(\frac{1}{2}z + \frac{1}{2}(x + 2td) \Big) \leq \frac{1}{2}F(z) + \frac{1}{2}\mu(2t) = \alpha,
\end{equation}
for a constant $\alpha$.  It follows that $F$ is bounded above by $\alpha$ on the ray $\{y + td, t \in \mathbb{R}^{+} \}$, from which it follows, by assumption, that $F$ is constant on the line $y + td$.
\end{proof}
We now prove Proposition \ref{firstlemma}. %\begin{proof}[Proof of Theorem \eqref{firstthm}]
Choosing $x_0 \in X$, we define the (nonempty) set
\begin{equation}
M := X \cap \{x \in \mathbb{R}^N, F(x) \leq F(x_0) \}.
\end{equation}
Obviously, the set $M$ is convex and closed.  By assumption, $F$ is bounded from below on $X$ and hence on $M$.  Therefore, if $M$ is bounded, then Weierstrass' Theorem yields the desired result.  
\\
\\
Thus, we may assume that $M$ is unbounded.  Then, the convexity of $M$ implies that $M$ contains a ray $r = \{z + td, t \geq 0\}$.   Denote by $r_1, r_2, ..., r_J$ a set of $J$ rays in $M$ corresponding to linearly independent vectors $d_1, ..., d_J$, so that any ray in $M$ can be expressed as a linear combination of the $r_1, ..., r_J$.   By definition of $M$ and by the assumption, $F$ is bounded on $M$, hence, $F$ is constant on each of the the lines $z_j + td_j$, according to Lemma \eqref{appendlemma1}.  From Lemma \eqref{appendlemma2}, it follows that $F$ is constant along each line $x + td_j$ for arbitrary $x \in \mathbb{R}^N$, from which we deduce that $F$ is constant along any line $x + td$ for arbitrary $d \in Y = \operatorname{span}\{d_1, ...., d_J \}$.  Thus, we project $X$ onto the subspace of $\mathbb{R}^N$ that is orthogonal to $Y$; call this subspace $\tilde{X}$.  
\\
\\
From the foregoing arguments, we have
\begin{equation}
\inf_{\tilde{X}} F(u) = \inf_X F(u)
\end{equation}
As $\tilde{X}$ is still a convex polyhedral set, and by construction $\tilde{M} = \tilde{X} \cap \{x \in \mathbb{R}^N \}$ contains no rays, Weierstrass' Theorem yields the desired result. 
%\end{proof}

\subsection{On uniform boundedness of $\| D^\dagger_h \|$}

\noindent The aim of the second part of the appendix is to prove the uniform bound $\| D^\dagger_h \| \leq 1/2$ eluded to in Section 3.1.  Again, $\| A \|$ denotes the spectral norm of the matrix $A$, and $D_h^\dagger: \mathbb{R}^{n-1} \rightarrow \mathbb R^n$ is the pseudo-inverse  of the discrete derivative matrix $D_h$ as given by \eqref{dermtrx}, with the identification $n = \lfloor 1/h \rfloor$.   From the expression for $D_h$, and the knowledge that $D_h D_h^\dagger = I$ is the identity operator and $D_h^\dagger D_h = (D_h^\dagger D_h)^*$ is self-adjoint, the $n \times (n-1)$ matrix $D_h^\dagger$ is identified as follows:
\begin{equation}
\label{dermtrxinverse}
D^{\dagger}_h = \frac{1}{n^2} \left ( \begin{array}{cccccc} -(n-1) & -(n-2) & -(n-3) & \dots & \dots & -1 \\
1 & -(n-2) & -(n-3) & \dots &\dots& -1 \\
1 & 2 & -(n-3) & \dots &\dots& -1 \\
\vdots & \vdots & \vdots & \vdots & \vdots & \vdots \\
1&2&3&\dots&\dots&n-1
\end{array} \right ).
\end{equation}
It is well-known that the spectral norm of an $m \times n$ matrix can be bounded by the more manageable entry-wise Frobenius norm, according to
\begin{equation}
\| A \| \leq \| A \|_F = \sqrt{\sum_{i=1}^m \sum_{j=1}^n | a_{i,j}|^2}.
\end{equation}
As such, we need only to bound the sum of the squares of the entries of $D^{\dagger}_h$.  The sum $S^1_n = \sum_{j=1}^{n-1} | d_{1,j}|^2$ over entries in the first row of $D^{\dagger}_h$ is given by $S^1_n = (n-1)(2n-1)/(6n^3)$, using the familiar formula $\sum_{j=1}^N j^2 = \frac{1}{6}N(N+1)(2N+1)$.  The analogous sum over entries in the $j^{th}$ row of $D^{\dagger}_h$ is seen inductively to satisfy $S^j_n = S^1_n - \frac{(j-1)}{n^2} +\frac{j(j-1)}{n^3}$.  The total sum $S_n = \sum_{j=1}^{n} S_n^j$ is then $S_n = \frac{1}{6} - \frac{1}{6n^2}$, and we arrive at the desired uniform bound:
\begin{equation}
\| D_h^\dagger \| \leq \sqrt{S_n} \leq \frac{1}{\sqrt{6}} < 1/2. 
\end{equation}  

\subsection{Proof of Proposition \ref{thmthrs}}

\noindent In order to help the reading of the current proof, as well as the proofs of Theorem \ref{mainth} and Theorem \ref{globalmin} in later appendices, we report in Table 1 the notation of the functions used in the proof of Proposition \ref{thmthrs} for the definition of $H_{(p,r)}$.
\begin{table}[h]
\begin{center}
\begin{tabular}{|l|l|}
\hline
$L_p(t,\lambda)$   & $=(t - \lambda)^2 + \min\{ |t|^p, r^p \}$ \\
\hline
$G_p(t, \lambda)  $ & $=(t - \lambda)^2 + |t|^p$ \\
\hline
$F_p(t) $  & $=t+\frac{p}{2}\sgn t |t|^{p-1}$, $p >1$ \\
\hline
$S_p(\lambda) $  & $=G_p(F^{-1}_p(\lambda),\lambda) = (F^{-1}_p(\lambda) - \lambda)^2 + | F^{-1}_p(\lambda)|^p$, $p>1$\\
\hline
$H_{(p,r)}(\lambda)  $ & $=\arg \min_{t \geq 0} L_p(t, \lambda)$ for general $\lambda \geq 0$, $p>1$ \\
\hline
& $=\arg \min_{0 \leq t \leq r} G_p(t, \lambda)= F_p^{-1}(\lambda)$ for $0\leq \lambda  \leq r$ \\
\hline
& $=\left\{ \begin{array}{ll}  F_p^{-1}(\lambda), & \textrm{if } G_p(F_p^{-1}(\lambda),\lambda) \leq r^p \nonumber \\
\lambda, & \textrm{else} \end{array} \right.$ for $\lambda > r$.\\
\hline

\end{tabular}
\caption{Notation of the functions involved in the definition of $H_{(p,r)}$ as  in the proof of Proposition \ref{thmthrs}.}
\end{center}
\end{table}

%\begin{proof}[Proof of Theorem \eqref{thmthrs}] 
Consider the functions 
\begin{equation}
L_p(t,\lambda) = (t - \lambda)^2 + \min\{ |t|^p, r^p \},
\end{equation}
and
\begin{equation}
G_p(t, \lambda) = (t - \lambda)^2 + |t|^p.
\label{G}
\end{equation}
The proof reduces to solving for
\begin{equation}
H_{(p,r)}(\lambda) = \arg \min_{t \in \mathbb R} L_p(t, \lambda)
\label{hmin}
\end{equation}
as a function of $\lambda \in \mathbb{R}$.  Since $L_p(t, \lambda) = L_p(-t, -\lambda)$, the function $H_{(p,r)}(\lambda)$ will be odd, and since also $H_{(p,r)}(0) = 0$, we can, without loss of generality, restrict the domain of interest to $\lambda > 0$.  On this domain, $H_{(p,r)}(\lambda) = \arg \min_{t \in \mathbb R} L_p(t, \lambda)$ is nonnegative, since $L_p(t, \lambda) \leq L_p(-t, \lambda)$ when $t \geq 0$ and $\lambda \geq 0$. Hence, we can restrict the minimization of $L_p(t,\lambda)$ to $t \geq 0$.
\\

\noindent It will be convenient to split the proof into two cases: $1 < p$ and $p = 1$.
\begin{enumerate}
\item We first analyze the case $1 < p$.  
\\
Note that 
\begin{eqnarray}
\arg \min_{t \geq r} L_p(t, \lambda) &=&  \arg \min_{t \geq r}  (t - \lambda)^2 \nonumber \\
&=& \max \{\lambda, r \},
\end{eqnarray}
so that the minimization  $\eqref{hmin}$ naturally splits into the following two cases:
\begin{enumerate}
\item If $\lambda \leq r$, the minimizer has to be searched in $[0,r]$, hence
\begin{equation}
H_{(p,r)}(\lambda) = \arg \min_{0 \leq t \leq r} G_p(t, \lambda) = F_p^{-1}(\lambda) \leq \lambda 
\label{lessthan}
\end{equation}
where $F_p^{-1}(\lambda)$ is the functional inverse of the increasing, and continuous function
\begin{equation}
F_p(t) = t + \frac{p}{2}\sgn{t}{|t|}^{p-1}.
\end{equation}
\item On the other hand, if $\lambda > r$, the minimizer has to be searched in $[0, \lambda]$, hence
\begin{eqnarray}
H_{(p,r)}(\lambda) &=& \left\{ \begin{array}{ll}  F_p^{-1}(\lambda), & \textrm{if } G_p(F_p^{-1}(\lambda),\lambda) \leq r^p \nonumber \\
\lambda, & \textrm{else} \end{array} \right. .
\label{h2}
\end{eqnarray}
\end{enumerate}
By implicit differentiation of the functional relation $F_p(F_p^{-1}(\lambda)) = \lambda$, it is clear that the functions $F_p^{-1}(\lambda)$ and $S_p(\lambda) := G_p(F_p^{-1}(\lambda),\lambda)$ are strictly increasing functions in $\lambda$.  
Indeed, we have the bounds
$$
0 < \frac{d}{d \lambda}F_p^{-1}(\lambda) = \left (F'_p(F_p^{-1}(\lambda)) \right)^{-1} = \left (1 + \frac{p(p-1)}{2} (F_p^{-1}(\lambda))^{p-2} \right )^{-1} \leq 1,
$$
and
\begin{eqnarray*}
\frac{d}{d \lambda} S_p(\lambda) &=& \frac{\partial}{\partial t} G_p(F_p^{-1}(\lambda),\lambda)\frac{d}{d \lambda}F_p^{-1}(\lambda) + \frac{\partial}{\partial \lambda}   G_p(F_p^{-1}(\lambda),\lambda)\\
&=& (2 (F_p^{-1}(\lambda) - \lambda) + p (F_p^{-1}(\lambda))^{p-1})\frac{d}{d \lambda}F_p^{-1}(\lambda) - 2 (  F_p^{-1}(\lambda) - \lambda)\\
&=& 2 \left (1 - \frac{d}{d \lambda}F_p^{-1}(\lambda)\right ) (\lambda - F_p^{-1}(\lambda)) + p \frac{d}{d \lambda}F_p^{-1}(\lambda) (F_p^{-1}(\lambda))^{p-1} \geq 0,
\end{eqnarray*}
since $0\leq \frac{d}{d \lambda}F_p^{-1}(\lambda) \leq 1$, and 
\begin{equation}
\label{boundinvF}
0\leq F_p^{-1}(\lambda) \leq \lambda.
\end{equation}
Also observe that $F_p^{-1}(r + \frac{p}{2}r^{p-1}) = r$, and $S_p(r + \frac{p}{2}r^{p-1}) = r^p + \frac{p^2}{4}r^{2p-2} > r^p$.
This leads us to immediately conclude that
\begin{itemize}
\item[(i)] If $\lambda \leq r $, then $H_{(p,r)}(\lambda) = F_p^{-1}(\lambda)$ (from $\eqref{lessthan}$).
\item[(ii)] If $\lambda \geq r +\frac{p}{2}r^{p-1}$, then $S_p(\lambda) = G_p(F_p^{-1}(\lambda),\lambda) > r^p$, so that $H_{(p,r)}(\lambda) = \lambda$.  
\item[(iii)] Since $S_p(r) < r^p$ while $S_p(r+\frac{p}{2} r^{p-1})) > r^p$, the intermediate value theorem implies that there exists a unique value $\lambda'(r,p)$  lying {\it strictly within} the interval $\big(r , r^{p-1}(\frac{p}{2} + r^{2-p})\big)$ at which 
 \begin{equation}
S_p(\lambda')  = r^p,
\label{equal}
 \end{equation} 
 and  
 \begin{eqnarray}
 H_{(p,r)}(\lambda) &=& \left\{ \begin{array}{ll}  F_p^{-1}(\lambda) & \lambda < \lambda'(r,p)  \\ \lambda & \lambda > \lambda'(r,p) \end{array} \right. .
 \end{eqnarray}
At $\lambda'$, $H_{(p,r)}(\lambda') = \arg \min_{t \geq 0} L_p(t, \lambda')$ is not uniquely defined and is realized at $F_p^{-1}(\lambda')$ and at $\lambda'$. In this case, we identify $H_{(p,r)} (\lambda') = F_p^{-1}(\lambda)$ for the sequel; as will be made clear, this will not cause problems in the ensuing analysis.  Finally, note that
 \item[(iv)] 
 At $\lambda'$, the function $H_{(p,r)}$ has a discontinuity $\delta(r,p) = \lambda' - H_{(p,r)} (\lambda')$ that is strictly positive, as long as $r > 0$.  Indeed,
on the one hand, we know that $\lambda'(r,p) > r$, on the other hand, $H_{(p,r)} (\lambda') < r$. This follows because $H_{(p,r)}(\lambda') = F_p^{-1}(\lambda')$, and
$$
(F_p^{-1}(\lambda'))^p < (F_p^{-1}(\lambda') - \lambda')^2 + |F_p^{-1}(\lambda')|^p = S_p(\lambda')= r^p.
$$
 \end{itemize}

 \item The analysis of the case $p = 1$ is left to the reader since it follows a similar argument as for $p > 1$.
 \end{enumerate}
%\end{proof}

\subsection{Proof of Theorem \ref{contract} }

%\begin{proof}.
We assume that the operator $T^*T : \ell_2(\mathcal I) \rightarrow \ell_2(\mathcal I)$ is nonnegative, so that its spectrum lies within an interval $[\delta, 1]$ with $\delta \geq 0$, and the operator $I - T^*T$ has norm $\| I - T^*T\| \leq 1 - \delta$.  In particular, if $T^*T$ is invertible, then the inequality $\delta > 0$ is strict, and so $\|I - T^*T\| \leq 1 - \delta < 1$.
\\
We wish to show that the map $\mathbb{U}_{\mathcal I_0}$ with component-wise action
\begin{eqnarray}
[\mathbb{U}_{\mathcal I_0} u]_i &=& \left\{
\begin{array}{ll} F_p^{-1}([(I - T^*T)u + T^*g]_i), & \textrm{if } i \in \mathcal I_0 \\  
\big( (I - T^*T)u + T^*g \big)_i, & \textrm{if } i \in \mathcal I_1
\end{array} \right.
\label{separate}
\end{eqnarray}
is a contraction.  To this end, let ${ v, v'}$ be arbitrary vectors in $\ell_2(\mathcal I)$.  
\begin{enumerate}
\item If the index $i \in \mathcal I_1$, then $$|[\mathbb{U}_{\mathcal I_0}(v)]_i - [\mathbb{U}_{\mathcal I_0}(v')]_i|  = |\big[ (I - T^*T)(v - v') \big]_i |;$$
\item  If the index $i \in \mathcal I_0$, then we split the analysis in two cases $p>1$ and $p=1$:
\begin{enumerate}
\item for $p > 1$, we have
\begin{eqnarray}
|[\mathbb{U}_{\mathcal I_0}(v)]_i - [\mathbb{U}_{\mathcal I_0}(v')]_i| &=& \left | F_p^{-1}([(I - T^*T)v + T^*g]_i) - F_p^{-1}([(I - T^*T)v' + T^*g]_i) \right | \nonumber \\
&=& \left | \frac{d}{d\lambda}F_p^{-1}(\xi) \big[ (I - T^*T)(v - v') \big]_i \right | \nonumber \\
&<& \left |\big[ (I - T^*T)(v - v') \big]_i \right |
\end{eqnarray}
where the second equality is an application of the mean value theorem, which is valid since $F_p^{-1}(\lambda)$ is differentiable.  The final inequality above follows from implicit differentiation of the relation  $$F_p^{-1}(F_p(t)) = t $$
and the observation that $|\frac{d}{dt}F_p(t)| > 1$ (see the proof of Proposition \ref{thmthrs});
\item for $p = 1$, by analyzing all cases, we get also that
\begin{eqnarray}
|[\mathbb{U}_{\mathcal I_0}(v)]_i - [\mathbb{U}_{\mathcal I_0}(v')]_i| 
&\leq& |\big[ (I - T^*T)(v - v') \big]_i |
\end{eqnarray}
\end{enumerate}
\end{enumerate}
Together, we have 
\begin{eqnarray}
\| \mathbb{U}_{\mathcal I_0} ({ v}) - \mathbb{U}_{\mathcal I_0}({ v'}) \|^2_{\ell_2(\mathcal I)} &=& \sum_{i \in \mathcal I} |[\mathbb{U}_{\mathcal I_0}(v)]_i - [\mathbb{U}_{\mathcal I_0}(v')]_i|^2 \nonumber \\
&\leq&  \sum_{i \in \mathcal I} |\big( (I - T^*T) v - v' \big)_i |^2 \nonumber \\
&=& \| (I - T^*T) { v - v'} \|^2_{\ell_2(\mathcal I)} \nonumber \\
&\leq& \| I - T^*T \|^2 \|{ v - v'} \|^2_{\ell_2(\mathcal I)} \nonumber \\ 
&\leq& (1 - \delta) \|{ v - v'} \|^2_{\ell_2(\mathcal I)}.
\label{nonexpand}
\end{eqnarray}
As $\mathbb{U}_{\mathcal I_0}$ is a contraction, we arrive at the stated result by application of the Banach Fixed Point Theorem.
% \end{proof}

\subsection{Proof of Lemma \ref{l4}}
 If $i \in \mathcal I_1$, then $\bar{u}_i = \bar{u}_i + \big[ T^*(g - T\bar{u}) \big]_i$, which is satisfied if and only if $\big[ T^*(g - T\bar{u}) \big]_i = 0$ as stated.   It remains to analyze the case $i \in \mathcal I_0$, and, again, we split the argument in the cases $p>1$ and $p=1$.  
 \begin{enumerate}
 \item First suppose $p > 1$. Using the notation $\bar{\lambda} = \bar{u}_i + \big[ T^*(g - T\bar{u}) \big]_i$, the fixed point characterization $\eqref{seperate2}$ translates to $$ F_p^{-1}(\bar{\lambda}) = \bar{u}_i. $$
But of course $\lambda = F_p(\bar{u}_i)$ is the unique value at which $F_p^{-1}(\lambda) = \bar{u}_i$, and so this implies that
\begin{equation}
\big[ T^*(g - T\bar{u}) \big]_i = F_p(\bar{u}_i) - \bar{u}_i,
\label{invert}
\end{equation}
and, by reversing operations, the relation $\eqref{invert}$ in turn implies the fixed point condition $\eqref{seperate2}$.
\item The case $p = 1$, which is similar, is left to the reader.
\end{enumerate}

\subsection{Proof of Theorem \ref{mintheorem}}
The proof will be much simplified by the following lemma which characterizes vectors such as $\bar{u}$ that satisfy the fixed point relations \eqref{wha2} or \eqref{wah2}:
\begin{lemma}
If $u$ and $v$ are such that
\begin{equation}
{\cal{J}}_r^{p,surr} (u + v, u) - \| v \|_{\ell_2(\mathcal I)}^2\geq {\cal{J}}_r^{p,surr} (u, u) = {\cal{J}}_r^{p} (u), 
\label{eq}
\end{equation}
then ${\cal{J}}^p_r (u + v) \geq  {\cal{J}}^p_r ({ u})$.
\label{lemmata}
\end{lemma}
\begin{proof} 
For any $u$ and $v$, the following holds because $\|L\| \leq 1$:
\begin{equation}
{\cal{J}}^p_r(u + v) = {\cal{J}}_r^{p,surr} (u + v,u) - \|Lv \|_{\ell_2(\mathcal I)}^2 \geq  {\cal{J}}_r^{p,surr} (u + v, u) - \|v \|_{\ell_2(\mathcal I)}^2.
\end{equation}
If in addition $u$ and $v$ satisfy $\eqref{eq}$, then the desired result is achieved by virtue of the equality ${\cal{J}}_r^{p,surr} (u, u) =  {\cal{J}}^p_r ({ u})$.
\end{proof}

Let us show now the proof of Theorem \ref{mintheorem}.
By Lemma $\ref{lemmata}$, it suffices to show that at a fixed point $\bar{u}$ defined by \eqref{wha2} or \eqref{wah2}, any perturbation $\delta h \in \ell_2(\mathcal I)$ with norm $\| \delta h \|_{\ell_2(\mathcal I)} \leq \min\{ [\lambda'(r,p) - r], [r - H_{(p,r)}(\lambda')] \}$ will satisfy
\begin{equation}
{\cal{J}}_r^{p,surr} (\bar{u} + \delta h, \bar{u}) -{\cal{J}}_r^{p,surr} (\bar{u}, \bar{u}) \geq \| \delta h \|_{\ell_2(\mathcal I)}^2.
\label{min}
\end{equation}
After expanding the left-hand-side above, the inequality $\eqref{min}$ is seen to be equivalent to 
\begin{equation}
2 \sum_{i \in \mathcal I} \delta h_i [T^*(T\bar{u} - g)]_i + \sum_{i \in \mathcal I} \Big[ \min\{|\bar{u}_i + \delta h_i|^p, r^p\} - \min \{|\bar{u}_i|^p, r^p \} \Big] \geq 0.
\label{lhs}
\end{equation}
At this point, it is convenient to consider the summation over $i \in \mathcal I_0$ and $i \in \mathcal I_1$ separately. 
\\
\noindent By Lemma $\ref{l3}$, the first summand above vanishes over $\mathcal I_1$ and 
\begin{enumerate}
\item  if $1 < p $, then $\sum_{i \in \mathcal I} \delta h_i [T^*(T\bar{u} - g)]_i = - \sum_{i \in \mathcal I_0} \delta h_i \sgn{u_i} \frac{p}{2} |u_i|^{p-1}$;  
\item if $p = 1$, then $\sum_{i \in \mathcal I} \delta h_i [T^*(T\bar{u} - g)]_i = - 1/2 \sum_{i \in \mathcal I^b_0} \delta h_i \sgn{u_i} + \sum_{i \in \mathcal I^a_0} \delta h_i [T^*(T\bar{u} - g)]_i$.
\end{enumerate}
With respect to the second summation, observe from Proposition $\ref{thmthrs}$ that for all $1 \leq p$,  $|\bar{u}_i| \geq \lambda'(r,p) > r $ for $i \in \mathcal I_1$, so that this summation vanishes over $\mathcal I_1$ for any perturbation $\delta h$ satisfying the component-wise inequality $| \delta h_i | \leq \lambda'(r,p) - r$.  Similarly, $|\bar{u}_i| \leq H_{(p,r)} (\lambda') < r$ for $i \in \mathcal I_0$, so that for any perturbation $\delta h$ satisfying component-wise $| \delta h_i | \leq \min\{ [\lambda'(r,p) - r], [r - H_{(p,r)}(\lambda')] \}$, we have that
\begin{equation}
 \sum_{i \in \mathcal I} \Big[ \min\{|\bar{u}_i + \delta h_i|^p, r^p\} - \min \{|\bar{u}_i|^p, r^p \} \Big]  =  \sum_{i \in \mathcal I_0} |\bar{u}_i + \delta h_i|^p - |\bar{u}_i|^p.
\end{equation}
The desired result follows if we can show that
 \begin{enumerate}
 \item $1 < p \leq 2$:  $\big[ |\bar{u}_i + \delta h_i|^p - |\bar{u}_i|^p - \delta h_i p [\sgn{u_i}] |u_i|^{p-1} \big] \geq 0$, for all $i \in \mathcal I_0$
\item $p = 1$: 
\begin{enumerate}
\item $ | \delta h_i + \bar{u}_i| - |\bar{u}_i| - \delta h_i[\sgn{u_i}] \big] \geq 0$ for all $i \in \mathcal I_0^b$, and
\item  $\delta h_i [T^*(T\bar{u} - g)]_i + | \delta h_i | \geq 0$, for all $i \in \mathcal I_0^a$.
\end{enumerate}
\end{enumerate}
The inequality in $2(b)$ follows directly from Lemma $\ref{l3}$; by symmetry, $1$ and $2(a)$ follow if, for any $u \geq 0$,
\begin{equation}
\min_{v \in \mathbb R} \left [ f(v) := |u + v|^p - u^p - pu^{p-1}v \right ]  = \min_{v \geq -u} (u + v)^p - u^p - pu^{p-1}v \geq 0.
\end{equation}
When $p = 1$, the right-hand-side is identically zero and the result holds.  When $1 < p \leq 2$, differentiating the right-hand-side gives that $f(v)$ has a local minimum at $v = 0$, at which $f(0) = 0$, and, at the endpoint, $f(-u) = (p - 1)u^{p-1} \geq 0.$   
%\end{proof}

\subsection{Proof of Theorem \ref{globalmin}}
Suppose that $u^*$ is a minimizer of the functional ${\cal{J}}^p_r$.   Consider the partition of the index set $\mathcal I$ into $\mathcal I_0 = \{i \in \mathcal I: |u^*_i| \leq r \}$  and $\mathcal I_1 = \{i \in \mathcal I: |u^*_i| > r \}$, and note that $| \mathcal I_1 | < \infty$, or else $| \mathcal J_r^p(u^*)|$ would not be finite.  As in the proof of Theorem \eqref{l5}, consider the unique decomposition $u^* = u^*_0 + u^*_1$ into a vector $u^*_0$ supported on $\mathcal I_0$ and another $u^*_1$ supported on $\mathcal I_1$.  Again, let $\mathcal P: u \rightarrow u_1$ and $\mathcal P^{\perp} = \mathcal I - \mathcal P: u \rightarrow u_0$ denote the orthogonal projections onto the subspaces  $\ell_2^{\mathcal I_1}(\mathcal I)$ and $\ell_2^{\mathcal I_0}(\mathcal I)$, respectively, and consider the operators $T_0 = T \mathcal P^{\perp}$ and $T_1 = T \mathcal P$.
\\
\\
By minimality of ${ u^*}$, if we fix ${ u^*_0}$, the vector ${ u^*_1}$ satisfies ${ u_1^*} = \arg \min_{z \in \ell_2^{\mathcal I_1}(\mathcal I)} {\cal{J}}^p_{r,1}(z)$, where 
\begin{equation}
 {\cal{J}}^p_{r,1} ({ z}) :=  \| T_1 { z} - { (g - T_0 u^*_0)} \|^2_{\ell_2(\mathcal J)} +  \sum_{i \in \mathcal I_1} \min \{ |z_i|^p, r^p \}.
\label{restrict1}
\end{equation}
Since all coefficients in ${ u^*_1}$ have absolute value $|(u^*_1)_i| > r$, the vector ${ u^*_1}$ also minimizes the functional
\begin{eqnarray}
 \| T_1{ z} - { (g - T_0 u^*_0)} \|^2_{\ell_2(\mathcal J)}, 
\label{realrestrict}
\end{eqnarray}
or, else, the vector ${ z^*}$ minimizing $\eqref{realrestrict}$ would satisfy ${\cal{J}}^p_{r,1} ({ z^*}) <  {\cal{J}}^p_{r,1} ({ u^*_1})$, contradicting the minimality of ${ u^*_1}$.  In fact, ${ u^*_1}$ must be the {\it unique} vector minimizing $\eqref{realrestrict}$.  For, if another vector ${ u'}$ also minimized $\eqref{realrestrict}$, then the operator $T_1$ would have a nontrivial null space containing the span of some nonzero vector ${ v}$, so that all vectors in the affine space $\{ u^*_1 + t v : t \in \mathbb R\}$ would be minimal solutions for $\eqref{realrestrict}$.  In this case, we would have also the freedom of choosing from this affine subspace a vector ${ u'}$ having one coefficient $u'_i$ satisfying $|u_i'| < r$.  But such a vector ${ u'}$ satisfies ${\cal{J}}^p_{r,1} ({ u'}) <  {\cal{J}}^p_{r,1} ({ u^*_1})$, contradicting the minimality of $u^*_1$.  
\\
\\
\noindent It follows that the operator $T_1$ must have trivial null space, and ${ u^*_1}$ is the unique minimal least squares solution to $\eqref{realrestrict}$, well-known to be explicitly given by
\begin{equation}
{ u^*_1} = \big(T_1^*T_1\big)^{-1}T_1^*{ (g - T_0 u^*_0)},
\label{relu1u0}
\end{equation}
so that $T_1{ u^*_1}$ is the unique orthogonal projection of ${ (g - T_0 u^*_0)}$ onto the range of $T_1$. Actually $\mathcal P_1 = T_1 ( T_1^* T_1)^{-1} T_1^*$ is the orthogonal projection onto the range of $T_1$, due to the non-triviality of the null space of $T_1$.  Therefore we have $T_1{ u^*_1} = {\cal P}_1({ g - T_0 u^*_0})$. It easily follows that
\begin{equation}
T^*_1 \big(T_1{ u^*_1}  - { (g - T_0 u^*_0)} \big) = 0,
\end{equation}
or, in other words, 
\begin{equation}
 \big[ T^*(g - T{ u}^*) \big]_i  = 0 \textrm{,  for all }  i \in \mathcal I_1.
\end{equation}
Now, on the other hand, by observing that any optimal variable ${ u_1}$ for fixed $u_0$ depends on $u_0$ via the relationship ${ u_1} = \big(T_1^*T_1\big)^{-1}T_1^*{ (g - T_0 u_0)}$, we easily infer that the vector ${ u^*_0}$ minimizes
%\begin{eqnarray}
%{\cal{J}}^p_{r,0} ({ v}) &:=&  \| T_0 { v} - (g - T_1 u^*_1) \|^2_{\ell_2(\mathcal J)} +  \sum_{i \in \mathcal I_0} \min \{ |v_i|^p, r^p \}. \nonumber \\
%&=& \| {\cal P}_1^{\perp} (T_0 { v} - { g}) \|^2_{\ell_2(\mathcal J)} +  \sum_{i \in \mathcal I_0} \min \{ |v_i|^p, r^p \}. 
%\label{restrict1}
%\end{eqnarray}
\begin{eqnarray}
{\cal{J}}^p_{r,0} ({ v})
&=& \| {\cal P}_1^{\perp} (T_0 { v} - { g}) \|^2_{\ell_2(\mathcal J)} +  \sum_{i \in \mathcal I_0} \min \{ |v_i|^p, r^p \}, 
\label{restrict1}
\end{eqnarray}
where ${\cal P}_1^{\perp}$ denotes the orthogonal projection operator onto the orthogonal complement of the range of $T_1$.
\\
\\
\noindent Consider the convex functional,
\begin{equation}
 {\cal{F}} ({ v}) :=  \|  {\cal P}_1^{\perp} (T_0 { v} - { g}) \|^2_{\ell_2(\mathcal J)} + \|v \|_{\ell_p^{\mathcal I_0}(\mathcal I)}^p,
\label{auxfunc}
\end{equation}
and note that ${\cal{J}}^p_{r,0} ({ u}) \leq  {\cal{F}} ({ u})$, while at the same time ${\cal{J}}_{r,0}^p({ u^*_0}) =  {\cal{F}} ({ u^*_0})$ by virtue of the fact that $|u^*_i| < r$. 
{ For $p>1$} it follows that ${ u^*_0}$ is also a minimizer of  ${\cal{F}} ({ u})$, and so satisfies the Euler-Lagrange equations $\cite{bashsh93}$,
\begin{equation}
\big(T_0^*{\cal P}_1^{\perp}(T_0{ u^*_0}  - { g}) \big) + \frac{p}{2}\sgn{{ u^*_0}}|{ u^*_0}|^{p-1} = 0,
\end{equation}   
which imply the fixed point conditions
\begin{equation}
 \big[ T^*(g - T{ u}^*) \big]_i = \frac{p}{2}\sgn{(u^*_0)_i}|(u_0^*)_i|^{p-1}  \textrm{,  for all } i \in \mathcal I_0.
\end{equation}
For $p=1$ one uses results from \cite{DDD} to conclude that
\begin{equation}
u_0^* = \mathbb S_{1/2} ( u_0^* + T_0^* \mathcal P_1^\perp( g - T_0 u_0^*)),
\label{fixptst}
\end{equation}
where $\mathbb S_\gamma$ is the so-called {\it soft-thresholding}, defined component-wise $\mathbb S_\gamma(v) = (S_\gamma (v_i))_{i \in \mathcal I}$, where
\begin{equation}
S_{\gamma}(\lambda) =
\left\{
\begin{array}{ll}
0, & |\lambda| \leq \gamma  \\
\lambda - \frac{\sgn \lambda}{2}, & |\lambda|  > \gamma.
\end{array} \right.
\label{softthrs}
\end{equation}
(Actually, \cite[Proposition 3.10]{DDD} only states that any fixed point of \eqref{fixptst} is a minimizer of \eqref{auxfunc}; nevertheless the converse also holds, see \cite[Remarks (1), pag. 2515]{fo07}.)
The fixed-point condition \eqref{fixptst} implies
\begin{equation}
  \left\{
\begin{array}{ll}
\big[ T^*(g - T {u}^*) \big]_i  \in [-1/2, 1/2] , & |{u}^*_i| \leq 1/2 \\  
\big[ T^*(g - T{u}^*) \big]_i  = 1/2\sgn{{u}^*_i}, & 1/2 < |{u}^*_i| \leq r. 
\end{array}
\right.
\end{equation}
It remains to verify that 
\begin{itemize}
\item $|u^*_i| \geq \lambda'(r,p)$, if $i \in \mathcal I_1$, and
\item $|u^*_i| \leq F_p^{-1}(\lambda'(r,p))$, { for $p>1$, and $|u^*_i| \leq r- 1/4$, for $p=1$,} if $i \in \mathcal I_0$.
\end{itemize}
{ We show these conditions for $p>1$ only, as the case $p=1$ is proved with an analogous argument.}
\begin{enumerate}
\item
We first show that $|u^*_i| \geq \lambda'(r,p)$ if $i \in \mathcal I_1$. From the first part of the proof, we know that at a minimizer $u^*$, the functional ${\cal J}^p_r(u^*)$ can be written as
\begin{equation}
 {\cal J}^p_r(u^*) = \| {\cal P}_1^{\perp} (T_0 u^*_0 - { g}) \|^2_{\ell_2(\mathcal K)} +  \|u_0^* \|_{\ell_p^{\mathcal I_0}(\mathcal I)}^p+ |\mathcal I_1|r^p
 \label{atmin}
\end{equation} 
Note that at this point we make  explicit use of the finite cardinality of $\mathcal I_1$.
Fix $i \in \mathcal I_1$ and any perturbation ${ h} = h_i e_i$, $h_i \in \mathbb R$, along the coordinate $i$ (here, $e_i$ is the $i^{th}$ vector of the canonical basis).  Consider the rank-one operator $t_i = T \mathcal P_i$, where we use $\mathcal P_i$ to denote the orthogonal projection onto the one-dimensional subspace spanned by $e_i$.  Observe that $|\ t_i u \| = | (u)_i | \| t_i \|$.    Since $t_i$ is orthogonal to the argument ${\cal P}_1^{\perp} (T_0 u^*_0 - { g})$ under the $\ell_2$ penalty in $\eqref{atmin}$, the minimality condition ${\cal J}^p(u^*) \leq {\cal J}^p(u^* + h)$ can be written as
\begin{eqnarray}
&& \| {\cal P}_1^{\perp} (T_0 u^*_0 - { g}) \|^2_{\ell_2(\mathcal K)}+ \|u_0^* \|_{\ell_p^{\mathcal I_0}(\mathcal I)}^p+ |\mathcal I_1|r^p \nonumber \\
&\leq&    \| {\cal P}_1^{\perp} (T_0 u^*_0 - { g}) \|^2_{\ell_2(\mathcal K)} + \|u_0^* \|_{\ell_p^{\mathcal I_0}(\mathcal I)}^p \nonumber \\
&& \phantom{XXXXX} + \|h_i t_i\|^2_{\ell_2(\mathcal I)} + \min\{r^p, |u^*_i + h_i|^p\} + r^p (|\mathcal I_1| -1) 
\end{eqnarray}
which is equivalent to the condition that
\begin{eqnarray}
r^p &\leq& \|h_i t_i \|^2_{\ell_2(\mathcal I)} + \min\{r^p, |u^*_i + h_i|^p \}
\label{up}
\end{eqnarray}
hold for all $h_i \in \mathbb{R}$.  Now, since $\|T\|<1$, it follows that $\|t_i\| \leq 1$, and $\eqref{up}$ implies that
\begin{eqnarray}
r^p &\leq& h_i^2 + \min\{r^p, |u^*_i + h_i|^p \}
\label{down}
\end{eqnarray}
holds for all $h_i \in \mathbb{R}$, or, after the change of variables $\alpha = u^*_i + h_i$, that
\begin{equation}
r^p \leq (\alpha - u^*_i)^2 + \min\{r^p, |\alpha|^p \} 
\label{alpha}
\end{equation}
holds for all $\alpha \in \mathbb{R}$.  In particular, the inequality $\eqref{alpha}$ must hold at the value $\alpha^*$ that minimizes the right-hand-side.  But we already know from Proposition $\ref{thmthrs}$ that such a minimizer $\alpha^*$ is of the form:
\begin{eqnarray}
\alpha^* &=& \left\{
\begin{array}{ll}
F_p^{-1}(u^*_i), &  |u^*_i| \leq \lambda'(r,p)  \\
u^*_i, & |u^*_i|  > \lambda'(r,p) \\  
\end{array}
\right.
\end{eqnarray}
Now, suppose $|u^*_i| < \lambda'(r,p).$ (We know that $|u^*_i| > r$, so then $r < |u^*_i| < \lambda'(r,p)$).    From the proof of Proposition $\ref{thmthrs}$ we know that the function $F_p^{-1}(\lambda)$ is increasing, so then $\alpha^* = F_p^{-1}(u^*_i) < F_p^{-1}(\lambda') < r$.  Since also $S_p$ is strictly increasing, it follows that $S_p(\alpha^*) < S_p(F_p^{-1}(\lambda')) \leq S_p(\lambda') = r^p$. In the last inequality we used \eqref{boundinvF}. (See also Table 1 for recalling the notations used here.) 
But this is a contradiction to the minimality condition, $\eqref{alpha}$, and so we must conclude that $|u^*_i| \geq \lambda'(r,p)$.  
\item We now show that $|u^*_i| \leq F_p^{-1}(\lambda'(r,p))$, if $|u^*_i| \leq r$.   Recall that for $i \in \mathcal I_0$, the coefficient $u^*_i$ satisfies the fixed point condition,
\begin{equation}
 \big[ T^*(g - T{ u}^*) \big]_i  = \frac{p}{2}\sgn{u^*_i}|u^*_i|^{p-1}.
 \label{fp0}
\end{equation}
Fix $i \in \mathcal I_0$, and consider as before any perturbation ${ h} = h_i e_i$ along the coordinate $i$, $h_i \in \mathbb R$.  Let $t_i$ be the rank-one operator as defined before.  Then, the minimality condition ${\cal J}^p_r(u^*) \leq {\cal J}^p_r(u^* + h)$ is easily seen to be equivalent to
\begin{eqnarray}
\|Tu^* - g\|^2_{\ell_2(\mathcal K)} + |u^*_i|^p &\leq&  \|Tu^* - g + h_i t_i\|^2_{\ell_2(\mathcal K)} \nonumber \\
 &&+ \min\{r^p, |u^*_i + h_i|^p\} \nonumber \\
 &=& \|Tu^* - g \|^2_{\ell_2(\mathcal K)} + \|h_i t_i \|^2_{\ell_2(\mathcal I)} + 2 h_i \langle t_i,Tu^* - g \rangle \nonumber \\
 &&+ \min\{r^p, |u^*_i + h_i|^p\} \nonumber \\
 &=& \|Tu^* - g \|^2_{\ell_2(\mathcal K)} + \|h_i t_i \|^2_{\ell_2(\mathcal I)} - 2 h_i\frac{p}{2}\sgn{u^*_i}|u^*_i|^{p-1}  \nonumber \\
 && + \min\{r^p, |u^*_i + h_i|^p\} 
 \label{string}
\end{eqnarray}
and the final equality follows directly from the fixed point condition $\eqref{fp0}$.  Now the chain of inequalities $\eqref{string}$ implies the minimality condition
\begin{eqnarray}
|u^*_i|^p &\leq&  \|h_i t_i \|^2_{\ell_2(\mathcal I)} - 2h_i\frac{p}{2}\sgn{u^*_i}|u^*_i|^{p-1} + \min\{r^p, |u^*_i + h_i|^p\}  \nonumber \\
 &\leq& h_i^2  - 2h_i\frac{p}{2}\sgn{u^*_i}|u^*_i|^{p-1} + \min\{r^p, |u^*_i + h_i|^p\},
\end{eqnarray}
or, again using the change of variables  $\alpha = u^*_i + h_i$, the inequality
\begin{equation}
|u_i^*|^p \leq (\alpha - u^*_i)^2 - 2(\alpha - u^*_i)\frac{p}{2}\sgn{u^*_i}|u^*_i|^{p-1} + \min\{r^p, |\alpha|^p\}.
\label{again}
\end{equation}
%or
%\begin{equation}
%(1 - p)|u_i^*|^p \leq (\alpha - u^*_i)^2 - p \alpha \sgn{u^*_i}|u^*_i|^{p-1} + \min\{r^p, |\alpha|^p\}
%\label{again2}
%\end{equation}
Again, the inequality $\eqref{again}$ should hold for all $\alpha$ by the minimality of $u^*$.   Minimizers $\alpha^*$ of the right-hand-side of $\eqref{again}$ also are minimizers of  
\begin{equation}
\big(\alpha - (u^*_i + \frac{p}{2}\sgn{u^*_i}|u^*_i|^{p-1})\big)^2 + \min\{r^p, |\alpha|^p\},
\end{equation}
which we know to have the form
\begin{eqnarray}
\alpha^* &=& \left\{
\begin{array}{ll}
F_p^{-1}(u^*_i + \frac{p}{2}\sgn{u^*_i}|u^*_i|^{p-1}), &  |u^*_i| + \frac{p}{2}|u^*_i|^{p-1} \leq \lambda'(r,p)  \\
u^*_i + \frac{p}{2}\sgn{u^*_i}|u^*_i|^{p-1}, & |u^*_i| + \frac{p}{2}|u^*_i|^{p-1}  > \lambda'(r,p) \\  
\end{array}
\right. .
\end{eqnarray}
But $u^*_i + \frac{p}{2}\sgn{u_i^*}|u_i^*|^{p-1} = F_p(u_i^*)$, so the above reduces to
\begin{eqnarray}
\alpha^* &=& \left\{
\begin{array}{ll}
u^*_i, &  F_p(u_i^*) \leq \lambda'(r,p)  \\
F_p(u_i^*), & F_p(u^*_i) > \lambda'(r,p) \\  
\end{array}
\right.
\end{eqnarray}
As before, the proof proceeds by contradiction.  Suppose that $F_p(u_i^*)  > \lambda'(r,p)$, so that $\alpha^* = F_p(u_i^*) > \lambda'(r,p)$ and $S_p(\alpha^*) > S_p(\lambda') = r^p$.   Note that, by recalling $F_p(u_i^*) = u_i^* - \frac{p}{2} \sgn(u_i^*) (u_i^*)^{p-1}$, we have
\begin{equation}
 S_p(\alpha^*) = (u_i^* - F_p(u_i^*))^2 + |u_i^*|^p = |u_i^*|^p + \frac{p^2}{4}|u_i^*|^{2p-2}.
\label{bad}
\end{equation}  f
Plugging $\alpha^*$ into the right-hand-side of $\eqref{again}$, noting that $\lambda'(r,p) > r$ so that $|\alpha^*| \geq r$, and rearranging, yields the inequality
\begin{equation}
|u^*_i|^p \leq r^p -\frac{p^2}{4}|u_i^*|^{2p-2} \mbox{ or }  S_p(\alpha^*) = |u_i^*|^p + \frac{p^2}{4}|u_i^*|^{2p-2} \leq r^p.
\end{equation}
But this contradicts the assumption that the expression in $\eqref{bad}$ be larger than $r^p$.
\end{enumerate}
%\end{proof}

\subsubsection*{Acknowledgments}
We would like to thank Ingrid Daubechies and Albert Cohen  for various
conversations on the topic of this paper.  Massimo Fornasier acknowledges the financial support provided by the START-Prize ``Sparse Approximation and Optimization in High Dimensions'' of theFonds zur F\"orderung der wissenschaftlichen Forschung (FWF, Austrian Science Foundation), and he thanks the Program in Applied and Computational Mathematics at Princeton University for its hospitality during the early preparation of this work. The results of the paper also contribute to the project WWTF Five senses-Call 2006, Mathematical Methods for Image Analysis and Processing in the Visual Arts. \\ 
Rachel Ward acknowledges the hospitality of the Johann Radon Institute for Computational and Applied Mathematics, Austrian Academy of Sciences, for hosting her during the late preparation of this work.  She also acknowledges the support of the National Science Foundation Graduate Research Fellowship.

\bibliography{freesparsefinrev}

\end{document}